\documentclass{elsart}
\usepackage{amssymb}
\usepackage{graphicx}
\usepackage{amsmath}
\usepackage{color}
\usepackage{epstopdf}
\numberwithin{equation}{section}

\begin{document}
\begin{frontmatter}
\title{A multigrid correction scheme for a new Steklov eigenvalue problem in inverse scattering}
\author[author1,author2]{Yu Zhang},~~\ead{zhang\_hello\_hi@126.com}
\author[author1]{Hai Bi},~~\ead{bihaimath@gznu.edu.cn}
\author[author1]{Yidu Yang\corauthref{1}}~~\ead{ydyang@gznu.edu.cn}
\corauth[1]{Corresponding author}
\address[author1]{School of Mathematical Sciences, Guizhou Normal University, Guiyang 550001, China}
\address[author2]{School of Mathematics $\&$ Statistics, Guizhou University of Finance and Economics, Guiyang 550001, China}
\begin{abstract}
 We propose a multigrid correction scheme to solve a new Steklov eigenvalue problem in inverse scattering. With this scheme, solving an eigenvalue problem in a fine finite element space is reduced to solve a series of boundary value problems in fine finite element spaces and a series of eigenvalue problems in the coarsest finite element space. And the coefficient matrices associated with those linear systems are constructed to be symmetric and positive definite. We prove error estimates of eigenvalues and eigenfunctions. Numerical results coincide in theoretical analysis and indicate our scheme is highly efficient in solving the eigenvalue problem.
\end{abstract}
\end{frontmatter}
\section{Introduction}
\indent Steklov eigenvalue problems have important physical background and occur in many applications (see e.g., \cite{ahn,bergmann,bermudez,canavati,cakoni,cao,conca}). There have been various numerical approximation methods on Steklov eigenvalue problems (see \cite{armentano,armentano2,alonso,andreev,bramble,garau,hanx,li,lin,russo,tang,xie,xie2} and the references cited therein). Recently, a new Steklov eigenvalue problem in inverse scattering has attracted the attention of researchers (see \cite{cakoni,liu,bi}). \\
\indent Xu and Zhou propose a two grid method based on inverse iteration for elliptic eigenvalue problems in \cite{xu}. Later, it's developed to multigrid method (e.g., see \cite{dai,lin,xie0,yang}), among which \cite{lin,xie0} establish a new type of multigrid scheme based on the multilevel correction. And it's successfully applied to selfadjoint Steklov eigenvalue problem \cite{hanx,xie}, convecttion-diffusion eigenvalue problem \cite{peng}, transmission eigenvalue problem \cite{han}, etc. In the above applications, the associated bilinear or sesquilinear forms are coercive. In this paper the eigenvalue problem is non-selfadjoint and the associated  sesquilinear form is not $H^1$-ellipitic, which is the main difference from those studied before and causes the difficulty of theoretical analysis. Cakoni et al. study the conforming finite element approximation of this problem in \cite{cakoni}. Liu et al. prove the error estimate of eigenvalues in \cite{liu} for the first time and they prove that the discrete Neumann-to-Dirichlet operator $T_h$ converges to the Neumann-to-Dirichlet operator $T$ in the sense of norm $\|\cdot\|_{0, \partial \Omega}$. Furthermore, Bi et al. prove the convergence in the sense of norm $\|\cdot\|_{-\frac{1}{2}, \partial \Omega}$ in \cite{bi}.\\
\indent Based on the above work, in this paper, we present a multigrid correction scheme for the Steklov eigenvalue problem in inverse scattering and prove the error estimates of eigenvalues and eigenfunctions.  Without the loss of accuracy, solving eigenvalue problem in a fine finite element space is replaced by solving a series of boundary value problems in a series of fine finite element spaces and a series of eigenvalue problems in the coarsest finite element space. And the coefficient matrix associated with boundary value problem is constructed to be symmetric and positive definite. Numerical results coincide in our theoretical analysis and indicate this method is highly efficient for solving the Steklov eigenvalue problem in inverse scattering.\\
\indent The basic theory of finite element methods in this paper can be referred to \cite{babuska1,boffi,brenner,chatelin,ciarlet,oden,shi}. Throughout this paper, $C$ denotes a generic positive constant independent of mesh diameters, which may not be the same at each occurrence. For simplicity, we use the symbol $a\lesssim b$ to mean that $a\leqslant Cb$.
\section{Preliminary}
\indent Consider the following Steklov eigenvalue problem:
\begin{eqnarray}\label{s2.1}
  ~~~~~~~~~~~~~~~~~~~\Delta u+k^2n(x)u&=&0~~~~~in~~ \Omega,\\\label{s2.2}
  ~~~~~~~~~~~~~~~~~~\frac{\partial u }{\partial \nu}+\lambda u&=&0~~~~~on~~ \partial\Omega,
\end{eqnarray}
where $\Omega\subset \mathbb{R}^2$ is a bounded polygonal domain with Lipschitz boundary $\partial\Omega$, $\frac{\partial}{\partial\nu}$ denotes the unit outward normal derivative on $\partial\Omega$, $k$ is the wavenumber and $n(x)$ is the index of refraction. We assume $n(x)$ is a bounded complex value function given by:
\begin{eqnarray}\nonumber
  ~~~~~~~~~~~~n(x)=n_1(x)+i\frac{n_2(x)}{k},
\end{eqnarray}
and $i=\sqrt{-1}$, $n_1(x)\geq \delta>0$ and $n_2(x)\geq0$ are bounded and piecewise smooth functions.\\
\indent
Let $(\cdot, \cdot)_0$, $a(\cdot, \cdot)$ and $b(\cdot, \cdot)$ be defined as follows:
\begin{eqnarray*}
~~~~~~~~~~&&(u, v)_0=\int\limits_{\Omega} u\overline{v},\\
 ~~~~~~~~~~&&a(u,v)=(\nabla u, \nabla\overline{v})_0-(k^2n(x)u, v)_0,\\
 ~~~~~~~~~~&&b(u, v)=\int\limits_{\partial\Omega} uv.
\end{eqnarray*}
 \indent Thanks to \cite{cakoni}, the corresponding weak formulation to problem (\ref{s2.1})-(\ref{s2.2}) is given by: Find $(\lambda, u)\in \mathbb{C}\times H^1(\Omega)$, $\|u\|_{0, \partial\Omega}=1$, such that
\begin{eqnarray}\label{s2.3}
  ~~~~~~~~~~~~&~&a(u, v)=-\lambda b(u, v),~~~~ \forall v\in H^1(\Omega).
\end{eqnarray}
\indent Let $\pi_h=\{K\}$ be a shape-regular decomposition of  $\Omega$ into triangles. $h_K$ denotes the diameter of $K$. Let $h = \max\limits_{K\in \pi_h}\{h_K\}$ be the diameter of $\pi_h$.
$$V_h=\{v\in H^1(\Omega): v|_K\in P_1, \forall~K\in \pi_h\},$$
where $P_1$ denotes the space of linear polynomials. $\partial V_h$ denotes the restriction of $V_h$ on $\partial\Omega$. \\
\indent The finite element approximation of (\ref{s2.3}) is to find $(\lambda_h, u_h)\in \mathbb{C}\times V_h$, $\|u_h\|_{0, \partial\Omega}=1$, such that
\begin{eqnarray}\label{s2.4}
  ~~~~~~~~~~~~&~&a(u_h, v)=-\lambda_h b(u_h, v),~~~~ \forall v\in V_h.
\end{eqnarray}
From Section 2 in \cite{bi} we know, for any $f\in H^{-\frac{1}{2}}(\partial \Omega)$  the operator $A: H^{-\frac{1}{2}}(\partial \Omega)\rightarrow H^1(\Omega)$ can be defined as
\begin{equation}\label{s2.5}
  a(Af, v)=b(f, v),~~~~\forall~~v\in H^1(\Omega).
\end{equation}
and the Neumann-to-Dirichlet operator $T: H^{-\frac{1}{2}}(\partial \Omega)\rightarrow H^{\frac{1}{2}}(\partial \Omega)$ by
\begin{equation*}
  Tf=Af|_{\partial \Omega}.
\end{equation*}
Similarly, define a discrete operator $A_h: H^{-\frac{1}{2}}(\partial \Omega)\rightarrow V_h$ as
\begin{equation}\label{s2.6}
  a(A_hf, v)=b(f, v),~~~~\forall~~v\in V_h,
\end{equation}
and the discrete Neumann-to-Dirichlet operator $T_h: H^{-\frac{1}{2}}(\partial \Omega)\rightarrow \partial V_h$ such that
\begin{equation*}
  T_hf=A_hf|_{\partial \Omega}.
\end{equation*}
Thus (\ref{s2.3}) and (\ref{s2.4}) has the following equivalent operator form respectively:
\begin{eqnarray}\label{s2.7}
  ~~~~~~~~~~~~~~~~~~~~~Au &=& -\frac{1}{\lambda} u,~~~~~ Tu=-\frac{1}{\lambda} u\\\label{s2.8}
  ~~~~~~~~~~~~~~~~~~~~~A_hu_h &=& -\frac{1}{\lambda_h} u_h,~~~ T_hu_h = -\frac{1}{\lambda_h} u_h.
\end{eqnarray}
One can define $\eta_0(h)$ as
\begin{equation}\label{s2.9}
  \eta_0(h)=\sup\limits_{f\in H^\frac{1}{2}(\partial\Omega), \|f\|_{\frac{1}{2}, \partial\Omega}=1}\inf\limits_{v\in V_h}\|Af-v\|_{1,\Omega}.
\end{equation}
Define $P_h$ be the finite element projection operator of $H^1(\Omega)$ onto $V_h$ satisfying
\begin{equation}\label{s2.10}
  a(w-P_hw, v)=0,~~~~\forall~w\in H^1(\Omega)~~ and ~~ v\in V_h.
\end{equation}
\indent Consider the dual problem of (\ref{s2.3}): Find $(\lambda^*, u^*)\in \mathbb{C}\times H^1(\Omega)$, $\|u^*\|_{0, \partial\Omega}=1$, such that
\begin{eqnarray}\label{s2.11}
  ~~~~~~~~~~~~&~&a(v, u^*)=-\overline{\lambda^*} b(v, u^*),~~~~ \forall v\in H^1(\Omega).
\end{eqnarray}
The primal and dual eigenvalues are connected via $\lambda=\overline{\lambda^*}$.\\
\indent The finite element approximation associated with (\ref{s2.11}) is given by: Find $(\lambda^*_h, u^*_h)\in \mathbb{C}\times V_h$, $\|u^*_h\|_{0, \partial\Omega}=1$, such that
\begin{eqnarray}\label{s2.12}
  ~~~~~~~~~~~~&~&a(v, u^*_h)=-\overline{\lambda^*_h} b(v, u^*_h),~~~~ \forall v\in V_h.
\end{eqnarray}
The primal and dual eigenvalues are connected via $\lambda_h=\overline{\lambda^*_h}$.\\
\indent Similarly, from source problems corresponding to (\ref{s2.11}) and (\ref{s2.12}) we can define the operators $A^*: H^{-\frac{1}{2}}(\partial \Omega)\rightarrow H^1(\Omega)$ and $A_h^*: H^{-\frac{1}{2}}(\partial \Omega)\rightarrow V_h$ such that
\begin{eqnarray}\label{s2.13}
   ~~~~~~~~~~~~~~ a(v, A^*f)&=&b(v, f),~~~~\forall~~v\in H^1(\Omega), \\ \label{s2.14}
   ~~~~~~~~~~~~~~ a(v, A^*_hf)&=& b(v, f) ~~~~\forall~~v\in V_h.
\end{eqnarray}
Analogously, Neumann-to-Dirichlet and discrete Neumann-to-Dirichlet operators can be defined  by $T^*: H^{-\frac{1}{2}}(\partial \Omega)\rightarrow H^{\frac{1}{2}}(\partial \Omega)$ and $T_h^*: H^{-\frac{1}{2}}(\partial \Omega)\rightarrow \partial V_h$.\\
Furthermore, one can define $\eta^*_0(h)$
\begin{equation}\label{s2.15}
  \eta^*_0(h)=\sup\limits_{f\in H^\frac{1}{2}(\partial\Omega), \|f\|_{\frac{1}{2}, \partial\Omega}=1}\inf\limits_{v\in V_h}\|A^*f-v\|_{1,\Omega}.
\end{equation}
Let $P^*_h: H^1(\Omega)\rightarrow V_h$  be the projection defined by
\begin{equation}\label{s2.16}
  a(v, w^*-P^*_hw^*)=0,~~~~\forall~w^*\in H^1(\Omega)~~ and ~~ v\in V_h.
\end{equation}
\indent There holds the following lemma for boundary value problem (\ref{s2.5}), which is needed in our theory analysis.\\
{\bf Lemma 2.1}.~~ If $f\in L^2(\partial \Omega)$, then $Af\in H^{1+\frac{\gamma}{2}}(\Omega)$ and
\begin{equation}\label{s2.17}
  \|Af\|_{1+\frac{\gamma}{2}}\leq C\|f\|_{0, \partial\Omega},
\end{equation}
if $f\in H^\frac{1}{2}(\partial \Omega)$, then $Af\in H^{1+\gamma}(\Omega)$ and
\begin{equation}\label{s2.18}
  \|Af\|_{1+\gamma}\leq C\|f\|_{\frac{1}{2}, \partial\Omega},
\end{equation}
where $\gamma=1$ when the largest inner angle $\theta$ of $\Omega$ satisfying $\theta<\pi$, and $\gamma<\frac{\pi}{\theta}$ which can  be arbitrarily close to $\frac{\pi}{\theta}$ when $\theta>\pi$.\\
{\bf Proof.}~~See \cite{dauge}. ~~~$\Box$\\
\indent The dual problem (\ref{s2.13}) has the same regularity as the corresponding source problem.\\
\indent Refer to Lemma 2.2 in \cite{bi}, we can prove the following lemma.\\
{\bf Lemma 2.2}.~~For $\forall w\in H^1(\Omega)$, the following estimates hold:
\begin{eqnarray}\label{s2.19}
  ~~~~~~~~~~~~&~& \eta_0(h)\rightarrow 0,~~~\eta^*_0(h)\rightarrow 0,~~(h\rightarrow 0),\\ \label{s2.20}
 ~~~~~~~~~~~~ &~&\|w-P_hw\|_{-\frac{1}{2}, \partial\Omega}\lesssim\eta^*_0(h)\|w-P_hw\|_{1, \Omega},\\\label{s2.21}
  ~~~~~~~~~~~~&~&\|w^*-P^*_hw^*\|_{-\frac{1}{2}, \partial\Omega}\lesssim\eta_0(h)\|w^*-P^*_hw^*\|_{1, \Omega}.
\end{eqnarray}
{\bf Proof.}~~(\ref{s2.19}) can be proved by (\ref{s2.18}) and the standard error estimate of interpolation. Next, we prove (\ref{s2.20}).\\
\indent According to (\ref{s2.13}) and (\ref{s2.10}) we deduce
\begin{eqnarray*}
  \|w-P_hw\|_{-\frac{1}{2}, \partial\Omega} &=& \sup\limits_{f\in H^\frac{1}{2}(\partial\Omega), \|f\|_{\frac{1}{2}, \partial\Omega}=1}|b(w-P_hw, f)|\\
     &=&\sup\limits_{f\in H^\frac{1}{2}(\partial\Omega), \|f\|_{\frac{1}{2}, \partial\Omega}=1}|a(w-P_hw, A^*f)|\\
    &=&\sup\limits_{f\in H^\frac{1}{2}(\partial\Omega),\|f\|_{\frac{1}{2}, \partial\Omega}=1}|a(w-P_hw, A^*f-v)|\\
    &\lesssim& \sup\limits_{f\in H^\frac{1}{2}(\partial\Omega), \|f\|_{\frac{1}{2}, \partial\Omega}=1 }\|w-P_hw\|_{1,\Omega}\|A^*f-v\|_{1,\Omega},~~~\forall~v\in V_h
\end{eqnarray*}
which combines with the definition of $\eta^*_0(h)$ to yield (\ref{s2.20}). Similarly, (\ref{s2.21}) can be proved.~~~~$\Box$\\
\indent We define the following sesquilinear form:
\begin{equation*}
 \tilde{a}(u, v)=(\nabla u, \nabla v)_0+(u, v)_0,
\end{equation*}
\indent According to the definition of $\tilde{a}$, (\ref{s2.3}) can be written as
\begin{equation}\label{s2.22}
  \tilde{a}(u, v)=-\lambda b(u, v)+((k^2n(x)+1)u, v)_0,~~~\forall v\in H^1(\Omega).
\end{equation}
\indent Due to $\tilde{a}(\cdot, \cdot)$ is $H^1$-elliptic, we can define Ritz projection $\tilde{P}_h, \tilde{P}^*_h:~H^1(\Omega)\rightarrow V_h$  which satisfies
\begin{eqnarray}\label{s2.23}
   ~~~~~~~~~~~~~~  \tilde{a}(w-\tilde{P}_hw, v)=0,~~~~\forall~w\in H^1(\Omega)~~ and ~~ v\in V_h,\\\label{s2.24}
    ~~~~~~~~~~~~~~ \tilde{a}(v, w^*-\tilde{P}^*_hw^*)=0,~~~~\forall~w\in H^1(\Omega)~~ and ~~ v\in V_h.
\end{eqnarray}
\indent For our later discussion, we consider auxiliary problems as follows:
Find $\psi, \psi^*\in H^1(\Omega)$  such that
\begin{eqnarray*}
  ~~~~~~~~~~~~~~~~~ \tilde{a}(\psi, v)&=&b(f, v), ~~~\forall~v\in H^1(\Omega),\\
   ~~~~~~~~~~~~~~~~~ \tilde{a}(v, \psi^*)&=&b(v, f), ~~~\forall~v\in H^1(\Omega).
\end{eqnarray*}
For any $f\in H^{-\frac{1}{2}}(\partial \Omega)$  the operator $A_1, A_1^* : H^{-\frac{1}{2}}(\partial \Omega)\rightarrow H^1(\Omega)$ satisfies respectively
\begin{eqnarray}\label{s2.25}
   ~~~~~~~~~~~~~~ \tilde{a}(A_1f, v)=b(f, v),~~~~\forall~~v\in H^1(\Omega),\\ \label{s2.26}
   ~~~~~~~~~~~~~~ \tilde{a}(v,A_1^*f)=b(v, f),~~~~\forall~~v\in H^1(\Omega).
\end{eqnarray}
One can define $\eta_1(h)$ and $\eta^*_1(h)$ as
\begin{eqnarray*}
  ~~~~~~~~~~~~~~\eta_1(h)=\sup\limits_{f\in H^\frac{1}{2}(\partial\Omega), \|f\|_{\frac{1}{2}, \partial\Omega}=1}\inf\limits_{v\in V_h}\|A_1f-v\|_{1,\Omega},\\
 ~~~~~~~~~~~~~~ \eta^*_1(h)=\sup\limits_{f\in H^\frac{1}{2}(\partial\Omega), \|f\|_{\frac{1}{2}, \partial\Omega}=1}\inf\limits_{v\in V_h}\|A_1^*f-v\|_{1,\Omega}.
\end{eqnarray*}
Next, we prove the following lemma  holds.\\
{\bf Lemma 2.3.} For $\forall w\in H^1(\Omega)$, the following estimates hold:
\begin{eqnarray}\label{s2.27}
 ~~~~~~~~~~~~ &~& \eta_1(h)\rightarrow 0, ~~~ \eta^*_1(h)\rightarrow 0,~~~(h\rightarrow 0), \\\label{s2.28}
 ~~~~~~~~~~~~ &~&\|w-\tilde{P}_hw\|_{-\frac{1}{2}, \partial\Omega}\lesssim\eta^*_1(h)\|w-\tilde{P}_hw\|_{1, \Omega},\\\label{s2.29}
  ~~~~~~~~~~~~&~&\|w^*-\tilde{P}^*_hw^*\|_{-\frac{1}{2}, \partial\Omega}\lesssim\eta_1(h)\|w^*-\tilde{P}^*_hw^*\|_{1, \Omega}.
\end{eqnarray}
{\bf Proof.}~~Using similar arguments to (\ref{s2.19}) we can prove (\ref{s2.27}).
(\ref{s2.28}) follows (\ref{s2.26}), (\ref{s2.23}) and  the definition of $\eta^*_1(h)$
Similarly, combining (\ref{s2.25}), (\ref{s2.24}) and the definition of $\eta_1(h)$ we deduce (\ref{s2.30})~~~~$\Box$\\
\indent For our later discussion, we consider auxiliary problems:
Find $\psi_f, \psi^*_f\in H^1(\Omega)$  such that
\begin{eqnarray*}
  ~~~~~~~~~~~~~~ \tilde{a}(\psi_f, v)&=&(f, v)_0, ~~~\forall~v\in H^1(\Omega),\\
   ~~~~~~~~~~~~~~ \tilde{a}(v, \psi^*_f)&=&(v, f)_0, ~~~\forall~v\in H^1(\Omega).
\end{eqnarray*}
For any $f\in L^2(\Omega)$  the operator $A_2, A_2^* : L^2(\Omega)\rightarrow H^1(\Omega)$ satisfies respectively
\begin{eqnarray}\label{s2.30}
 ~~~~~~~~~~~~~~\tilde{a}(A_2f, v)=(f, v)_0,~~~~\forall~~v\in H^1(\Omega),\\ \label{s2.31}
 ~~~~~~~~~~~~~~\tilde{a}(v,A_2^*f)=(v, f)_0,~~~~\forall~~v\in H^1(\Omega).
\end{eqnarray}
One can define $\eta_2(h)$ and $\eta^*_2(h)$ as
\begin{eqnarray*}
  ~~~~~~~~~~~~~~\eta_2(h)=\sup\limits_{f\in L^2(\Omega), \|f\|_{0, \Omega}=1}\inf\limits_{v\in V_h}\|A_2f-v\|_{1,\Omega},\\
 ~~~~~~~~~~~~~~  \eta^*_2(h)=\sup\limits_{f\in L^2(\Omega), \|f\|_{0, \Omega}=1}\inf\limits_{v\in V_h}\|A_2^*f-v\|_{1,\Omega}.
\end{eqnarray*}
Then the following lemma holds.\\
{\bf Lemma 2.4.} For $\forall w\in H^1(\Omega)$, the following estimates hold:
\begin{eqnarray}\label{s2.32}
 ~~~~~~~~~~~~&~& \eta_2(h)\rightarrow 0,~~~\eta^*_2(h)\rightarrow 0,~~~(h\rightarrow 0),\\ \label{s2.33}
 ~~~~~~~~~~~~ &~&\|w-\tilde{P}_hw\|_{0, \Omega}\lesssim\eta^*_2(h)\|w-\tilde{P}_hw\|_{1, \Omega},\\\label{s2.34}
  ~~~~~~~~~~~~&~&\|w^*-\tilde{P}^*_hw^*\|_{0, \Omega}\lesssim\eta_2(h)\|w^*-\tilde{P}^*_hw^*\|_{1, \Omega}.
\end{eqnarray}
{\bf Proof.}~~Using similar arguments to (\ref{s2.19}) we can prove (\ref{s2.32}).
(\ref{s2.33}) follows (\ref{s2.31}), (\ref{s2.23}) and  the definition of $\eta^*_1(h)$
Similarly, combining (\ref{s2.30}), (\ref{s2.24}) and the definition of $\eta_1(h)$ we deduce (\ref{s2.34})~~~~$\Box$\\
\indent  Liu et al. prove $\|T_h-T\|_{0, \partial \Omega}\rightarrow 0 (h\rightarrow 0)$ in \cite{liu}. Furthermore, Bi et al. prove $\|T_h-T\|_{-\frac{1}{2}, \partial \Omega}\rightarrow 0 (h\rightarrow 0)$ in \cite{bi}. Actually it easy to prove that $\|A_h-A\|_{1, \Omega}\rightarrow 0 (h\rightarrow 0)$.\\
\indent Let $\lambda$ be the
$i$th eigenvalue of (\ref{s2.3}) with the ascent $\alpha$ and the algebraic multiplicity $q$. Then there are $q$ eigenvalues $\lambda_{j,h}(j=i,i+1,\cdots, i+q-1)$ of (\ref{s2.4}) converging to $\lambda$.  Let $M(\lambda)$ be the space spanned by all generalized eigenfunctions of (\ref{s2.3}) corresponding to the eigenvalue $\lambda$. Let $M_h(\lambda)$ be the space spanned
by all generalized eigenfunctions of (\ref{s2.4}) corresponding to the eigenvalues $\lambda_{j,h}(j=i,i+1,\cdots, i+q-1)$. As
for the dual problems (\ref{s2.11}) and (\ref{s2.12}), the definitions of $M^*(\lambda^*)$ and $M^*_h(\lambda^*)$ are made similarly to $M(\lambda)$ and $M_h(\lambda)$, respectively. Then we define
\begin{eqnarray}\label{s2.35}
  ~~~~~~~~~~~~~~~~~~\delta_h(\lambda)&=&\sup\limits_{w\in M(\lambda), \|w\|_{0, \partial\Omega}=1}\inf\limits_{v\in V_h}\|w-v\|_{1,\Omega},\\\label{s2.36}
  ~~~~~~~~~~~~~~~~~~\delta^*_h(\lambda^*)&=&\sup\limits_{w^*\in M^*(\lambda^*),  \|w^*\|_{0, \partial\Omega}=1}\inf\limits_{v\in V_h}\|w^*-v\|_{1,\Omega}.
\end{eqnarray}
\indent Thanks to Lemmas 2.2, 2.4 and 2.6 in \cite{bi} and the spectral approximation theory (see \cite{babuska1}) we get the following conclusions.\\
{\bf Lemma 2.5}.~~Assume $u_h$ is the eigenfunction approximation of (\ref{s2.4}), then there exists an eigenfunction of (\ref{s2.3}) $u\in M(\lambda)$ corresponding to $\lambda$ such that
\begin{eqnarray}\label{s2.37}
~~~~~~~~~~~~~~~~~&~&\|u-u_h\|_{-\frac{1}{2}, \partial\Omega}\lesssim \big\{\eta^*_0(h)\delta_h(\lambda)\big\}^\frac{1}{\alpha},\\\label{s2.38}
~~~~~~~~~~~~~~&~&\|u-u_h\|_{1,\Omega}\lesssim \delta_h(\lambda)+\big\{\eta^*_0(h)\delta_h(\lambda)\big\}^\frac{1}{\alpha},\\\label{s2.39}
~~~~~~~~~~~~~~~~~&~&\|u-u_h\|_{0, \Omega}\lesssim\big\{\eta^*_0(h)\delta_h(\lambda)\big\}^\frac{1}{\alpha},
\end{eqnarray}
and
\begin{equation}\label{s2.40}
|\lambda-\lambda_h|\lesssim \big\{\delta_h(\lambda)\delta^*_h(\lambda^*)\big\}^\frac{1}{\alpha}.
\end{equation}
{\bf Lemma 2.6}.~~Assume $u^*_h$ is the eigenfunction approximation of (\ref{s2.12}), then there exists an eigenfunction of (\ref{s2.11}) $u^*\in M^*(\lambda^*)$ corresponding to $\lambda^*$ such that
\begin{eqnarray}\label{s2.41}
~~~~~~~~~~~~~~~~~~&~&\|u^*-u^*_h\|_{-\frac{1}{2}, \partial\Omega}\lesssim \big\{\eta_0(h)\delta^*_h(\lambda^*)\big\}^\frac{1}{\alpha},\\\label{s2.42}
~~~~~~~~~~~~~~~~~&~&\|u^*-u^*_h\|_{1,\Omega}\lesssim \delta^*_h(\lambda^*)+\big\{\eta_0(h)\delta^*_h(\lambda^*)\big\}^\frac{1}{\alpha},\\\label{s2.43}
~~~~~~~~~~~~~~~~~&~&~\|u^*-u^*_h\|_{0, \Omega}\lesssim \big\{\eta_0(h)\delta^*_h(\lambda^*)\big\}^\frac{1}{\alpha},
\end{eqnarray}
and
\begin{equation}\label{s2.44}
|\lambda^*-\lambda^*_h|\lesssim\big\{\delta_h(\lambda)\delta^*_h(\lambda^*)\big\}^\frac{1}{\alpha}.
\end{equation}
\section{One correction step}
\indent In this section, based on the work in \cite{lin,xie,xie2}, we establish Algorithm 3.1 (One Correction Step). Firstly, initial mesh is given by $\pi_H=\pi_{h_1}$ with the mesh size $H=h_1$. We define a sequence of triangulation $\pi_{h_{l+1}}$ with the mesh size $h_{l+1}$, which is produced by refining $\pi_{h_{l}}$ in the regular way. And
\begin{equation*}
  h_{l+1}\approx \frac{1}{\xi}h_{l},
\end{equation*}
where $\xi$ is an integer and always is 2 in our numerical experiments.\\
Based on the sequence of meshes, we define conforming linear finite element spaces as follows:
\begin{equation*}
  V_H=V_{h_1}\subset V_{h_2}\subset\cdots \subset V_{h_{n}}\subset H^1(\Omega),
\end{equation*}
then
\begin{equation}\label{s3.1}
  \delta_{h_{l+1}}(\lambda_j)\approx \frac{1}{\xi}\delta_{h_l}(\lambda_j),~~~\delta^*_{h_{l+1}}(\lambda^*_j)\approx \frac{1}{\xi}\delta^*_{h_l}(\lambda^*_j).
\end{equation}
\indent  Assume we have obtained the eigenpair approximations of (\ref{s2.3}) and (\ref{s2.11}) $(\lambda^c_{j, h_l}, u^c_{j, h_l})\in \mathbb{C}\times V_{h_l}$ with $\|u^c_{j, h_l}\|_{0, \partial\Omega}=1$  and  $(\lambda^{c*}_{j, h_l}, u^{c*}_{j, h_l})\in \mathbb{C}\times V_{h_l}$ with $\|u^{c*}_{j, h_l}\|_{0, \partial\Omega}=1$ for $j=i, i+1, \cdots, i+q-1$, respectively. Now, we give one correction step. \\
{\bf Algorithm 3.1.}~(One Correction Step)\\
{\bf Step 1}.~For $j=i,...,i+q-1$ , solve the following boundary problems:~~Find $\tilde{u}_{j, h_{l+1}}, \tilde{u}^*_{j, h_{l+1}}\in V_{h_{l+1}}$ such that
\begin{eqnarray}\label{s3.2}
\tilde{a}(\tilde{u}_{j, h_{l+1}}, v)&=&-\lambda^c_{j, h_{l}}b(u^c_{j, h_{l}}, v)+((k^2n(x)+1)u^c_{j, h_{l}}, v)_0,~~\forall v\in V_{ h_{l+1}},\\\label{s3.3}
\tilde{a}(v, \tilde{u}^*_{j, h_{l+1}})&=&-\overline{\lambda^{c*}_{j, h_{l}}}b(v, u^{c*}_{j, h_{l}})+(v, (k^2n(x)+1)u^{c*}_{j, h_{l}})_0,~~\forall v\in V_{ h_{l+1}}.
\end{eqnarray}
{\bf Step 2}.~Define a new finite element space:
\begin{equation*}
  V_{H,h_{l+1}}=V_H\bigoplus span \{\tilde{u}_{i, h_{l+1}},\cdots,\tilde{u}_{i+q-1, h_{l+1}}, \tilde{u}^*_{i, h_{l+1}},\cdots, \tilde{u}^*_{i+q-1, h_{l+1}}\},
\end{equation*}
and solve the following Steklov eigenvalue problem:\\
Find $(\lambda^c_{j, h_{l+1}}, u^c_{j, h_{l+1}}), (\lambda^{c*}_{j, h_{l+1}}, u^{c*}_{j, h_{l+1}}) \in \mathbb{C}\times V_{H, h_{l+1}}$ such that
\begin{eqnarray}\label{s3.4}
 ~~~~~~~~~~~~~ a(u^c_{j, h_{l+1}}, v)&=&-\lambda^c_{j, h_{l+1}}b(u^c_{j, h_{l+1}}, v),~~~\forall v\in V_{H, h_{l+1}}, \\\label{s3.5}
 ~~~~~~~~~~~~~ a(v, u^{c*}_{j, h_{l+1}})&=&-\overline{\lambda^{c*}_{j, h_{l+1}}}b(v, u^{c*}_{j, h_{l+1}}),~~~\forall v\in V_{H, h_{l+1}}.
\end{eqnarray}
\indent Output $\{\lambda^c_{j, h_{l+1}}\}^{i+q-1}_{j=i}$ and output $\{u^c_{j, h_{l+1}}\}^{i+q-1}_{j=i}\subset M_{h_{l+1}}(\lambda_{i})$ with $\|u^c_{j,h_{l+1}}\|_{0,\partial\Omega}=1$ and  $\{u^{c*}_{j, h_{l+1}}\}^{i+q-1}_{j=i}\subset M_{h_{l+1}}^{*}(\lambda_{i}^{*})$ with $\|u^{c*}_{j,h_{l+1}}\|_{0,\partial\Omega}=1$. We denote the two steps of Algorithm 3.1 by
\begin{eqnarray*}
  ~~~~~~~~~~~~~~ ~~&~&\{\lambda^c_{j, h_{l+1}},\lambda^{c*}_{j, h_{l+1}}, u^c_{j, h_{l+1}}, u^{c*}_{j, h_{l+1}}\}^{i+q-1}_{j=i}\\
&~&:=Correction(V_H, \{\lambda^c_{j, h_{l}},\lambda^{c*}_{j, h_{l}}, u^c_{j, h_{l}}, u^{c*}_{j, h_{l}}\}^{i+q-1}_{j=i}, V_{h_{l+1}}).
\end{eqnarray*}
 Generalized eigenfunction spaces $M_{h_{l+1}}(\lambda_{i})$ and $M^*_{h_{l+1}}(\lambda^*_{i})$ are defined as follows:
\begin{eqnarray*}
  ~~~~~~~~~~~M_{h_{l+1}}(\lambda_{i})&=&span\{u^c_{i, h_{l+1}},u^c_{i+1, h_{l+1}},\cdots, u^c_{i+q-1, h_{l+1}}\}, \\
  ~~~~~~~~~~~ M^*_{h_{l+1}}(\lambda^*_{i})&=&span\{u^{c*}_{i, h_{l+1}},u^{c*}_{i+1, h_{l+1}},\cdots, u^{c*}_{i+q-1, h_{l+1}}\}.
\end{eqnarray*}
\indent Let \\
$$\eta(H)=\max\{\eta_0(H), \eta_1(H), \eta_2(H)\},\eta^*(H)=\max\{\eta^*_0(H), \eta^*_1(H), \eta^*_2(H)\}.$$
 \indent We need to use the following assumption in order to make the error analysis.\\
\indent{(A0)} Suppose that there are $\{\tilde{u}_{j,h_{l}}\}^{i+q-1}_{j=i}\subset
M_{h_{l}}(\lambda_{i})$ with $\|\tilde{u}_{j,h_{l}}\|_{0, \partial\Omega}=1$,
and $\{\tilde{u}^{*}_{s,h_{l}}\}_{s=i}^{i+q-1}\subset
M^{*}_{h_{l}}(\lambda^{*}_{i})$ with $\|\tilde{u}^{*}_{s,h_{l}}\|_{0, \partial\Omega}=1$ such that
\begin{equation*}
  |b(u^c_{j,h_{l}},\tilde{u}^{*}_{s,h_{l}})|+
|b(\tilde{u}_{j,h_{l}},u^{c*}_{s,h_{l}})|\leq C(\eta(H)+\eta^{*}(H)),
\end{equation*}
where $j,s=i,\cdots,i+q-1, j\not=s$, and
$|b(u^c_{j,h_{l}},\tilde{u}^{*}_{j,h_{l}})|+|b(\tilde{u}_{j,h_{l}},u^{c*}_{j,h_{l}})| (j=i, i+1, \cdots, i+q-1)$ has a positive lower bound
uniformly with respect to $h_{l}$.\\
\indent In practical computing, we can use Arnoldi algorithm to solve the dual problem
 (\ref{s2.11}) and obtain $\{\tilde{u}^*_{s,h_{l}}\}^{i+q-1}_i$ and meanwhile MATLAB has provided the solvers ¡°sptarn¡± and ¡°eigs¡± to implement Arnoldi algorithm; we can also use the two sided Arnoldi algorithm in \cite{cullum} to compute both left and right eigenvectors of (\ref{s2.3}) at the same time, and obtain $\{u^c_{j,h_{l}}\}^{i+q-1}_i$, $\{\tilde{u}^*_{s,h_{l}}\}^{i+q-1}_i$.\\

{\bf Theorem 3.1}.~Assume (A0) holds and the accent $\alpha=1$. And there exist two numbers $\varepsilon_{h_{l}}(\lambda_{i})$,
 $\varepsilon^{*}_{h_{l}}(\lambda^{*}_{i})$ such that the given eigenpairs
 $\{\lambda^c_{j,h_{l}},u^c_{j,h_{l}}\}^{i+q-1}_{j=i}$ and $\{\lambda^{c*}_{j,h_{l}},u^{c*}_{j,h_{l}}\}_{j=i}^{i+q-1}$
 in Algorithm 3.1 have the following error estimates:
\begin{eqnarray}\label{s3.6}
  ~~~~~~~~~~~~~~~\|u^c_{j,h_l}-u_j\|_{1,\Omega} &\lesssim& \varepsilon_{h_l}(\lambda_i),\\\label{s3.7}
   ~~~~~~~~~~~~~~~\|u^{c*}_{j,h_l}-u^*_j\|_{1,\Omega} &\lesssim& \varepsilon^*_{h_l}(\lambda^*_i),\\\label{s3.8}
  ~~~~~~~~~~~~~~~\|u^c_{j,h_l}-u_j\|_{-\frac{1}{2}, \partial\Omega}&\lesssim&
  \eta^*(H)\varepsilon_{h_l}(\lambda_i),\\\label{s3.9}
    ~~~~~~~~~~~~~~~\|u^{c*}_{j,h_l}-u^*_j\|_{-\frac{1}{2}, \partial\Omega}&\lesssim&
    \eta(H)\varepsilon^*_{h_l}(\lambda^*_i),\\\label{s3.10}
  ~~~~~~~~~~~~~~~\|u^c_{j,h_l}-u_j\|_{0, \Omega}&\lesssim&
   \eta^*(H)\varepsilon_{h_l}(\lambda_i),\\\label{s3.11}
  ~~~~~~~~~~~~~~~\|u^{c*}_{j,h_l}-u^*_j\|_{0, \Omega}&\lesssim&
    \eta(H)\varepsilon^*_{h_l}(\lambda^*_i),\\\label{s3.12}
  ~~~~~~~~~~~~~~~|\lambda^c_{j,h_l}-\lambda_j|&\lesssim&
 \varepsilon_{h_l}(\lambda_i)\varepsilon^*_{h_l}(\lambda^*_i).
\end{eqnarray}
After implementing one correction step, the resultant approximation $\{\lambda^c_{j,h_{l+1}}, u^c_{j,h_{l+1}}\}^{i+q-1}_ {j=i}$ and $\{\lambda^{c*}_{j,h_{l+1}}, u^{c*}_{j,h_{l+1}}\}^{i+q-1}_{j=i}$ has error estimates as follows:
\begin{eqnarray}\label{s3.13}
 ~~~~~~~~~~~~~~~\|u^c_{j,h_{l+1}}-u_j\|_{1,\Omega} &\lesssim& \varepsilon_{h_{l+1}}(\lambda_i),\\\label{s3.14}
 ~~~~~~~~~~~~~~~\|u^{c*}_{j,h_{l+1}}-u^*_j\|_{1,\Omega} &\lesssim& \varepsilon^*_{h_{l+1}}(\lambda^*_i), \\\label{s3.15}
 ~~~~~~~~~~~~~~~\|u^c_{j,h_{l+1}}-u_j\|_{-\frac{1}{2}, \partial\Omega}&\lesssim&
 \eta^*(H)\varepsilon_{h_{l+1}}(\lambda_i),\\\label{s3.16}
 ~~~~~~~~~~~~~~~\|u^{c*}_{j,h_{l+1}}-u^*_j\|_{-\frac{1}{2}, \partial\Omega}&\lesssim&
\eta(H)\varepsilon^*_{h_{l+1}}(\lambda^*_i),\\\label{s3.17}
 ~~~~~~~~~~~~~~~\|u^c_{j,h_{l+1}}-u_j\|_{0, \Omega}&\lesssim&
 \eta^*(H)\varepsilon_{h_{l+1}}(\lambda_i),\\\label{s3.18}
 ~~~~~~~~~~~~~~~\|u^{c*}_{j,h_{l+1}}-u^*_j\|_{0, \Omega}&\lesssim&
 \eta(H)\varepsilon^*_{h_{l+1}}(\lambda^*_i),\\\label{s3.19}
 ~~~~~~~~~~~~~~~|\lambda^c_{j,h_{l+1}}-\lambda_j|&\lesssim&
 \varepsilon_{h_{l+1}}(\lambda_i)\varepsilon^*_{h_{l+1}}(\lambda^*_i),
\end{eqnarray}
where $\varepsilon_{h_{l+1}}(\lambda_j)=\eta^*(H)\varepsilon_{h_l}(\lambda_i)+\varepsilon_{h_l}(\lambda_i)\varepsilon^*_{h_l}(\lambda^*_i)+\delta_{h_{l+1}}(\lambda_i)$ and $\varepsilon^*_{h_{l+1}}(\lambda^*_j)=\eta(H)\varepsilon^*_{h_l}(\lambda^*_i)+\varepsilon_{h_l}(\lambda_i)\varepsilon^*_{h_l}(\lambda^*_i)+\delta^*_{h_{l+1}}(\lambda^*_i)$. \\
{\bf Proof.}~~Since $\{u_j\}_{j=i}^{i+q-1}$ is a
basis of $M(\lambda_{i})$, for any $w\in M(\lambda_{i})$, $\|w\|_{0,\partial\Omega}=1$ we denote
\begin{equation}\label{s3.20}
w=\sum^{i+q-1}\limits_{j=i}\gamma_{j}u_{j}.
\end{equation}
Then
\begin{equation*}
b(w,u_{s}^{*})=\sum\limits_{j=i}^{i+q-1}\gamma_{j}b(u_{j},u_{s}^{*}).
\end{equation*}
Hence
\begin{equation}\label{s3.21}
\gamma_{s}=\frac{1}{b(u_{s},u_{s}^{*})}\big\{b(w,u_{s}^{*})-\sum\limits_{j\not=s,j=i}^{i+q-1}\gamma_{j}b(u_{j},u_{s}^{*})\big\},~~~s=i,i+1,\cdots,i+q-1.
\end{equation}
From (\ref{s3.8}), (\ref{s3.9}) and assumption (A0) we have
\begin{eqnarray}\nonumber
 |b(u_{j},u_{s}^{*})|&\leq& |b(u_{j}-u^c_{j,h_{l}},u_{s}^{*})|+|b(u^c_{j,h_{l}},u_{s}^{*}-u_{s,h_{l}}^{c*})|+|b(u^c_{j,h_{l}},u_{s,h_{l}}^{c*})|\\\label{s3.22}
  &\lesssim&
  \eta^*(H)\varepsilon_{h_l}(\lambda_i)+\eta(H)\varepsilon^*_{h_l}(\lambda^*_i)+\eta(H)+\eta^{*}(H).
\end{eqnarray}
Due to assumption (A0), we know $|b(u_{s}, u^*_{s})|$ has a positive lower bound uniformly with respect to $h_{l}$, i.e., there exists $\delta>0$ such that
\begin{equation*}
|b(u_{s},u_{s}^{*})|\geq \delta,~~~s=i,i+1,\cdots,i+q-1.
\end{equation*}
Substituting the above relation and (\ref{s3.22}) into (\ref{s3.21}), we conclude
\begin{eqnarray}\nonumber
|\gamma_{s}|&\leq& \frac{1}{|b(u_{s},u_{s}^{*})|}\big\{|b(w,u_{s}^{*})|+\sum\limits_{j\not=s,j=i}^{i+q-1}|\gamma_{j}||b(u_{j},u_{s}^{*}|)\big\}\\\nonumber
&\lesssim& \frac{1}{\delta}\big\{1+\sum\limits_{j\not=s,j=i}^{i+q-1}|\gamma_{j}|(\eta^*(H)\varepsilon_{h_l}(\lambda_i)+\eta(H)\varepsilon^*_{h_l}(\lambda^*_i)+\eta(H)+\eta^{*}(H))\big\}.
\end{eqnarray}
From which it follows that
\begin{equation}\nonumber
\sum\limits_{s=i}^{i+q-1}|\gamma_{s}|\lesssim
\frac{q}{\delta}\big\{1+\sum\limits_{j=i}^{i+q-1}|\gamma_{j}|\eta^*(H)\varepsilon_{h_l}(\lambda_i)+\eta(H)\varepsilon^*_{h_l}(\lambda^*_i)+\eta(H)+\eta^{*}(H))\big\}.
\end{equation}
The above inequality shows that
\begin{equation}\label{s3.23}
\sum\limits_{j=i}^{i+q-1}|\gamma_{j}|\lesssim 1.
\end{equation}
We set $\alpha_{j}:=\frac{\lambda_{i}}{\lambda^c_{j,h_{l}}}(j=i,...,i+q-1)$. According to the ellipticity of $\tilde{a}(\cdot, \cdot)$, (\ref{s2.23}), (\ref{s3.2}) and (\ref{s2.22}), we have
\begin{eqnarray}\nonumber
&&\|\alpha_{j}\tilde{u}_{j,h_{l+1}}-\tilde{P}_{h_{l+1}}u_{j}\|^{2}_{1,\Omega}\\\nonumber
&&~\lesssim |\tilde{a}(\alpha_{j}\tilde{u}_{j,h_{l+1}}-\tilde{P}_{h_{l+1}}u_{j},\alpha_{j}\tilde{u}_{j,h_{l+1}}-\tilde{P}_{h_{l+1}}u_{j})|\\\nonumber
&&~\lesssim |\tilde{a}(\alpha_{j}\tilde{u}_{j,h_{l+1}},\alpha_{j}\tilde{u}_{j,h_{l+1}}-\tilde{P}_{h_{l+1}}u_{j})-\tilde{a}(u_j,\alpha_{j}\tilde{u}_{j,h_{l+1}}-\tilde{P}_{h_{l+1}}u_{j})|\\\nonumber
&&~\lesssim |\lambda_jb(u^c_{j,h_{l}}-u_j, \alpha_{j}\tilde{u}_{j,h_{l+1}}-\tilde{P}_{h_{l+1}}u_{j})|\\\nonumber
&&~~~~+|\alpha_{j}-1|((k^2n(x)+1)u^c_{j, h_l},\alpha_{j}\tilde{u}_{j,h_{l+1}}-\tilde{P}_{h_{l+1}}u_{j})_0|\\\nonumber
&&~~~~+|((k^2n(x)+1)(u^c_{j, h_l}-u_j), \alpha_{j}\tilde{u}_{j,h_{l+1}}-\tilde{P}_{h_{l+1}}u_{j})_0|\\\nonumber
&&~\lesssim (\|u^c_{j, h_{l}}-u_{j}\|_{{-\frac{1}{2}},\partial\Omega}+|\lambda^c_{j, h_{l}}-\lambda_j|+\|u^c_{j, h_{l}}-u_{j}\|_{{0},\Omega})\|\alpha_{j}\tilde{u}_{j,h_{l+1}}-\tilde{P}_{h_{l+1}}u_{j}\|_{1,\Omega}.\\\label{s3.24}
\end{eqnarray}
Substituting (\ref{s3.8}), (\ref{s3.12}) and (\ref{s3.10}) into (\ref{s3.24}) we obtain
\begin{equation}\nonumber
\|\alpha_{j}\tilde{u}_{j,h_{l+1}}-\tilde{P}_{h_{l+1}}u_{j}\|_{1,\Omega}\lesssim
\eta^*(H)\varepsilon_{h_l}(\lambda_i)+\varepsilon_{h_l}(\lambda_i)\varepsilon^*_{h_l}(\lambda^*_i).
\end{equation}
Based on the above inequality and the error estimate of finite element projection $\|u_j-\tilde{P}_{h_{l+1}}u_{j}\|_{1,\Omega}\lesssim \delta_{h_{l+1}}(\lambda_i)$, we deduce
\begin{eqnarray}\label{s3.25}
&&\|\alpha_{j}\tilde{u}_{j,h_{l+1}}-u_j\|_{1,\Omega}\lesssim
\eta^*(H)\varepsilon_{h_l}(\lambda_i)+\varepsilon_{h_l}(\lambda_i)\varepsilon^*_{h_l}(\lambda^*_i)+\delta_{h_{l+1}}(\lambda_i).
\end{eqnarray}
Using (\ref{s3.23}) and (\ref{s3.25}), we deduce
\begin{eqnarray}\nonumber
~~~~~~~~~~~&&\sup\limits_{w\in
M(\lambda_{i}),\|w\|_{0, \partial\Omega}=1}\inf\limits_{v\in
V_{H,h_{l+1}}}\|w-v\|_{1, \Omega}\\\nonumber
&&~~~~~~\lesssim \sup\limits_{w\in
M(\lambda_{i}), \|w\|_{0,\partial\Omega}=1}\|\sum\limits_{j=i}^{i+q-1}\gamma_{j}u_{j}
-\sum\limits_{j=i}^{i+q-1}\gamma_{j}\alpha_{j}\tilde{u}_{j,h_{l+1}}\|_{1, \Omega}\\\nonumber
&&~~~~~~\lesssim \sup\limits_{\gamma_{j}}\|\sum\limits^{i+q-1}_{j=i}\gamma_{j}(u_{j}-\alpha_{j}\tilde{u}_{j,h_{l+1}})\|_{1, \Omega}\\\nonumber
&&~~~~~~\lesssim \mathrm{max}\big\{\{\|\alpha_{j}\tilde{u}_{j,h_{l+1}}-u_{j}\|_{1,\Omega}\}^{i+q-1}_{j=i}\big\}\\\nonumber
&&~~~~~~\lesssim
\eta^*(H)\varepsilon_{h_l}(\lambda_i)+\varepsilon_{h_l}(\lambda_i)\varepsilon^*_{h_l}(\lambda^*_i)+\delta_{h_{l+1}}(\lambda_i)\\\label{s3.26}
&&~~~~~~\lesssim \varepsilon_{h_{l+1}},
\end{eqnarray}
where $\varepsilon_{h_{l+1}}:=\eta^*(H)\varepsilon_{h_l}(\lambda_i)+\varepsilon_{h_l}(\lambda_i)\varepsilon^*_{h_l}(\lambda^*_i)+\delta_{h_{l+1}}(\lambda_i)$.\\
Now we estimate the error of $ u^c_{j, h_{l+1}}$.\\
Define $\tilde{\eta}^*_1(H)$ by
\begin{equation*}
  \tilde{\eta}^*_1(H)=\sup\limits_{f\in H^\frac{1}{2}(\partial\Omega), \|f\|_{\frac{1}{2}, \partial\Omega}=1}\inf\limits_{v\in V_{H,h_{l+1}}}\|A_1^*f-v\|_{1,\Omega}.
\end{equation*}
 It's easy to know $\tilde{\eta}^*_1(H) \lesssim \eta^*_1(H)$.\\
 From spectral approximation theory (see \cite{babuska1}), (\ref{s2.28}) and (\ref{s3.26}) we have
\begin{eqnarray}\nonumber
&&\|u^c_{j, h_{l+1}}-u_j\|_{-\frac{1}{2}, \partial\Omega}\lesssim\sup\limits_{w\in
M(\lambda_{i}),\|w\|_{0, \partial\Omega}=1}\inf\limits_{v\in
V_{H,h_{l+1}}}\|w-v\|_{-\frac{1}{2}, \partial\Omega}\\ \nonumber
&&~~~\lesssim \tilde{\eta}^*_1(H)\sup\limits_{w\in
M(\lambda_{i}),\|w\|_{0, \partial\Omega}=1}\inf\limits_{v\in
V_{H,h_{l+1}}}\|w-v\|_{1,\Omega}\\\label{s3.27}
&&~~~\lesssim \eta^*_1(H)\varepsilon_{h_{l+1}}.
\end{eqnarray}
Using similar proof to (2.26) in \cite{bi} and (\ref{s3.26}), we have
\begin{eqnarray}\label{s3.28}
&&\|u^c_{j, h_{l+1}}-u_j\|_{1, \Omega}\lesssim \sup\limits_{w\in
M(\lambda_{i}),\|w\|_{0, \partial\Omega}=1}\inf\limits_{v\in
V_{H,h_{l+1}}}\|w-v\|_{1, \Omega}\lesssim \varepsilon_{h_{l+1}}.
\end{eqnarray}
By similar proof to (2.27) in \cite{bi}, (\ref{s2.33}) and (\ref{s3.26}), we obtain
\begin{equation}\label{s3.29}
\|u^c_{j, h_{l+1}}-u_j\|_{0, \Omega}
\lesssim \tilde{\eta}^*_2(H)\sup\limits_{w\in
M(\lambda_{i}),\|w\|_{0, \partial\Omega}=1}\inf\limits_{v\in
V_{H,h_{l+1}}}\|w-v\|_{1,\Omega}\lesssim \eta^*_2(H)\varepsilon_{h_{l+1}},
\end{equation}
where $\tilde{\eta}^*_2(H)=\sup\limits_{f\in L^2(\Omega), \|f\|_{0, \Omega}=1}\inf\limits_{v\in V_{H,h_{l+1}}}\|A_2^*f-v\|_{1,\Omega} \lesssim \eta^*_2(H)$.\\
Let $\eta^*(H):=\max\{\eta^*_0(H), \eta^*_1(H), \eta^*_2(H)\}$. (\ref{s3.13}), (\ref{s3.15}) and (\ref{s3.17}) follows from (\ref{s3.28}), (\ref{s3.27}) and (\ref{s3.29}), respectively.\\
\indent Analogously, conclusions (\ref{s3.14}), (\ref{s3.16}) and (\ref{s3.18}) hold.
By assumption (A0), (\ref{s3.13}) and (\ref{s3.14}) we obtain (\ref{s3.19}) is valid.~~~~$\Box$
\section{Multigrid Correction Scheme for the Steklov eigenvalue problem}
\indent In this section, we use the
correction step in the above section to establish a multigrid scheme for (\ref{s2.4}).\\
{\bf Algorithm 4.1.}~(Multigrid Correction Scheme)\\
{\bf Step 1}.~Construct a sequence of nested finite element spaces $V_H=V_{h_{1}},~V_{h_{2}},...,~V_{h_{n}}$ such that (\ref{s3.1}) holds. \\
{\bf Step 2}.~For $j=i, i+1, \cdots, i+q-1$, Solve the Steklov eigenvalue problems as follows:\\
find $(\lambda_{j,H}, u_{j,H})\in \mathbb{C}\times V_H$ such that $\|u_{j,H}\|_{0, \partial\Omega}=1$ and
\begin{equation*}
  a(u_{j,H}, v)=-\lambda_{j,H}b(u_{j,H}, v),~~~\forall v\in V_H,
\end{equation*}
find $(\lambda^*_{j,H}, u^*_{j,H})\in \mathbb{C}\times V_H$ such that $\|u^*_{j,H}\|_{0, \partial\Omega}=1$ and
\begin{equation*}
  a(u^*_{j,H}, v)=-\overline{\lambda^*_{j,H}}b(u^*_{j,H}, v),~~~\forall v\in V_H,
\end{equation*}
$\lambda^c_{j,h_1}=\lambda_{j,H}$, $u^c_{j,h_1}=u_{j,H}$, $\lambda^{c*}_{j,h_1}=\lambda^*_{j,H}$, $u^{c*}_{j,h_1}=u^*_{j,H}$.\\
{\bf Step 3}.~For $l=1, 2, \cdots, n-1$,
 \begin{eqnarray*}
 ~~~~~~~~~~&~& \{\lambda^c_{j, h_{l+1}},\lambda^{c*}_{j, h_{l+1}}, u^c_{j, h_{l+1}}, u^{c*}_{j, h_{l+1}}\}^{i+q-1}_{j=i}\\
 &~& =Correction(V_H, \{\lambda^c_{j, h_{l}},\lambda^{c*}_{j, h_{l}}, u^c_{j, h_{l}}, u^{c*}_{j, h_{l}}\}^{i+q-1}_{j=i}, V_{h_{l+1}})
\end{eqnarray*}
End.\\
We obtain $q$ eigenpair approximations $\{\lambda^c_{j, h_n}, u^c_{j, h_n}\}^{i+q-1}_{j=i}$, $\{\lambda^{c*}_{j, h_n}, u^{c*}_{j, h_n}\}^{i+q-1}_{j=i}\in \mathbb{C}\times V_{H, h_n}$ and $\lambda^c_{j, h_n}=\overline{\lambda^{c*}_{j, h_n}}$.\\
{\bf Theorem 4.1.}~~Suppose the conditions of Theorem 3.1 hold. Let the numerical eigenpairs $(\lambda^c_{j, h_n}, u^c_{j, h_n})$, $(\lambda^{c*}_{j, h_n}, u^{c*}_{j, h_n})(j=i, i+1, \cdots, i+q-1)$ be obtained by Algorithm 4.1. Then there exist eigenpairs $(\lambda_{j}, u_{j})$, $(\lambda^*_{j}, u^*_{j})$ such that the following estimates
hold
\begin{eqnarray}\label{s4.1}
 ~~~~~~~~~~~~~~~&&\|u^c_{j, h_{n}}-u_j\|_{1, \Omega}\lesssim \delta_{h_n}(\lambda_i),\\\label{s4.2}
~~~~~~~~~~~~~~~&&\|u^{c*}_{j, h_{n}}-u^*_j\|_{1, \Omega}\lesssim \delta^*_{h_n}(\lambda_i),\\\label{s4.3}
 ~~~~~~~~~~~~~~~&&\|u^c_{j, h_{n}}-u_j\|_{-\frac{1}{2}, \partial\Omega}\lesssim \eta^*(H)\delta_{h_n}(\lambda_i),\\\label{s4.4}
  ~~~~~~~~~~~~~~~&&\|u^{c*}_{j, h_{n}}-u^*_j\|_{-\frac{1}{2}, \partial\Omega}\lesssim \eta(H)\delta^*_{h_n}(\lambda^*_i),\\\label{s4.5}
 ~~~~~~~~~~~~~~~&&\|u^c_{j, h_{n}}-u_j\|_{0, \Omega}\lesssim \eta^*(H)\delta_{h_n}(\lambda_i),\\\label{s4.6}
 ~~~~~~~~~~~~~~~&&\|u^{c*}_{j, h_{n}}-u^*_j\|_{0, \Omega}\lesssim \eta(H)\delta^*_{h_n}(\lambda^*_i),\\\label{s4.7}
 ~~~~~~~~~~~~~~~&&|\lambda^c_{j, h_n}-\lambda_j|\lesssim \delta_{h_n}(\lambda_i)\delta^*_{h_n}(\lambda^*_i).
\end{eqnarray}
{\bf Proof.}~~ According to step 2 of Algorithm 4.1, Lemmas 2.5 and 2.6 we know
\begin{eqnarray*}
~~~~~~~~~~~~~~~~~\|u^c_{j,h_1}-u_j\|_{-\frac{1}{2}, \partial\Omega}&\lesssim&
\eta^*(H)\delta_{h_1}(\lambda_i),\\
~~~~~~~~~~~~~~~~~\|u^{c*}_{j,h_1}-u^*_j\|_{-\frac{1}{2}, \partial\Omega}&\lesssim&
\eta(H)\delta^*_{h_1}(\lambda^*_i),\\
~~~~~~~~~~~~~~\|u^c_{j,h_1}-u_j\|_{1,\Omega}&\lesssim&\delta_{h_1}(\lambda_i),\\
~~~~~~~~~~~~~~\|u^{c*}_{j,h_1}-u^*_j\|_{1,\Omega}&\lesssim& \delta^*_{h_1}(\lambda^*_i),\\
~~~~~~~~~~~~~~~~\|u^c_{j,h_1}-u_j\|_{0, \Omega}&\lesssim&\eta^*(H)\delta_{h_1}(\lambda_i),\\
~~~~~~~~~~~~~~~~~\|u^{c*}_{j,h_1}-u^*_j\|_{0, \Omega}&\lesssim&\eta(H)\delta^*_{h_1}(\lambda^*_i),\\
~~~~~~~~~~~~~~~~~|\lambda^c_{j,h_1}-\lambda_j|&\lesssim& \delta_{h_1}(\lambda_i)\delta^*_{h_1}(\lambda^*_i).
\end{eqnarray*}
Let  $\varepsilon_{h_{1}}(\lambda_{i}):=\delta_{h_{1}}(\lambda_{i})$. Noting $\{\varepsilon^*_{h_{m}}(\lambda^*_{i})\lesssim\eta^{*}(H)\}^n_{m=1}$ and using recursion we have
\begin{eqnarray*}
&&\varepsilon_{h_{n}}(\lambda_{i})=\eta^{*}(H)\varepsilon_{h_{n-1}}(\lambda_{i})+\varepsilon_{h_{n-1}}(\lambda_{i})\varepsilon^*_{h_{n-1}}(\lambda^*_{i})+\delta_{h_{n}}(\lambda_{i})\\
&~&~~~~~~~~~~~\lesssim
\eta^{*}(H)\varepsilon_{h_{n-1}}(\lambda_{i})+\delta_{h_{n}}(\lambda_{i})\\
&~&~~~~~~~~~~~\lesssim (\eta^{*}(H))^2\varepsilon_{h_{n-2}}(\lambda_{i})+\eta^{*}(H)\varepsilon_{h_{n-1}}(\lambda_{i})+\delta_{h_{n}}(\lambda_{i})\\
&~&~~~~~~~~~~~\lesssim \sum^n\limits_{l=1}(\eta^{*}(H))^{n-l}\delta_{h_{l}}(\lambda_{i})\\
&~&~~~~~~~~~~~\lesssim \sum^n\limits_{l=1}(\eta^{*}(H))^{n-l}\xi^{n-l}\delta_{h_{n}}(\lambda_{i})\\
&~&~~~~~~~~~~~\lesssim \frac{1}{1-\eta^{*}(H)\xi}\delta_{h_{n}}(\lambda_{i})\\
&~&~~~~~~~~~~~\lesssim \delta_{h_{n}}(\lambda_{i}).
\end{eqnarray*}
\indent Analogously, let $\varepsilon^{*}_{h_{1}}(\lambda^{*}_{i}):=\delta^{*}_{h_{1}}(\lambda^{*}_{i})$, we can prove $\varepsilon^*_{h_{n}}(\lambda^*_{i})\lesssim\delta^*_{h_{n}}(\lambda^*_{i})$.\\
\indent Using Theorem 3.1 we can obtain Theorem 4.1.~~~$\Box$
\section{Numerical experiments}
\indent In this section, in order to validate our theoretical results, the Multigrid Scheme (Algorithm 4.1) is applied to solve (\ref{s2.1})-(\ref{s2.2}) on three different domains (the square $(-\frac{\sqrt{2}}{2}, \frac{\sqrt{2}}{2})^2$, the L-shaped domain $(-1, 1)^2\setminus([0, 1)\times (-1, 0])$ and the square with a slit $(-\frac{\sqrt{2}}{2}, \frac{\sqrt{2}}{2})^2\setminus\{0\leq x\leq\frac{\sqrt{2}}{2}, y=0\}$). In computation, we select the index of refraction $n(x)=4$ or $n(x)=4+4i$.  For comparison, using linear element, we also solve the problem by the direct method. The discrete eigenvalue problems are solved in MATLAB 2016b on an Lenovo ideaPad PC with 1.8GHZ CPU and 8GB RAM. Our program is compiled under the package of iFEM \cite{chen}.  Since the exact eigenvalues are not known, we use the most accurate approximations in tables as reference eigenvalues. For convenience and simplicity, the
following notations are introduced in tables and figures.
\begin{description}
\item[]$h$:~~~~The diameter of meshes.
\item[]$\lambda_{j, h}$:~~The $j$th eigenvalue obtained by direct method on $\pi_h$.
\item[]$\lambda^c_{j,h}$:~~The $j$th eigenvalue obtained by Algorithm 4.1.
\item[]-- :~~~The calculation cannot proceed since the computer runs out of memory.
\end{description}

\indent Figs 1-3 depict the error curves for four eigenvalue approximations on each domain. On the figures, the closer the slope of one error curve is to -1, the closer the convergence order of corresponding eigenvalue approximation is to the optimal convergence order $O(h^2)$. According to the regularity theory, when ascent of $\lambda$ equals 1, we know the convergence order of the eigenvalue approximation $\lambda_{j,h}$ is $O(h^2)$ on the square. Not all convergence order of the eigenvalue approximation $\lambda_{j,h}$ can reach $O(h^2)$ on the L-shaped domain and the square with a slit. Figs 1-3 indicate that the numerical results are coincided in the theory.
 \begin{figure}[h!]
 \centering
\includegraphics[width=6.5cm,height=5cm]{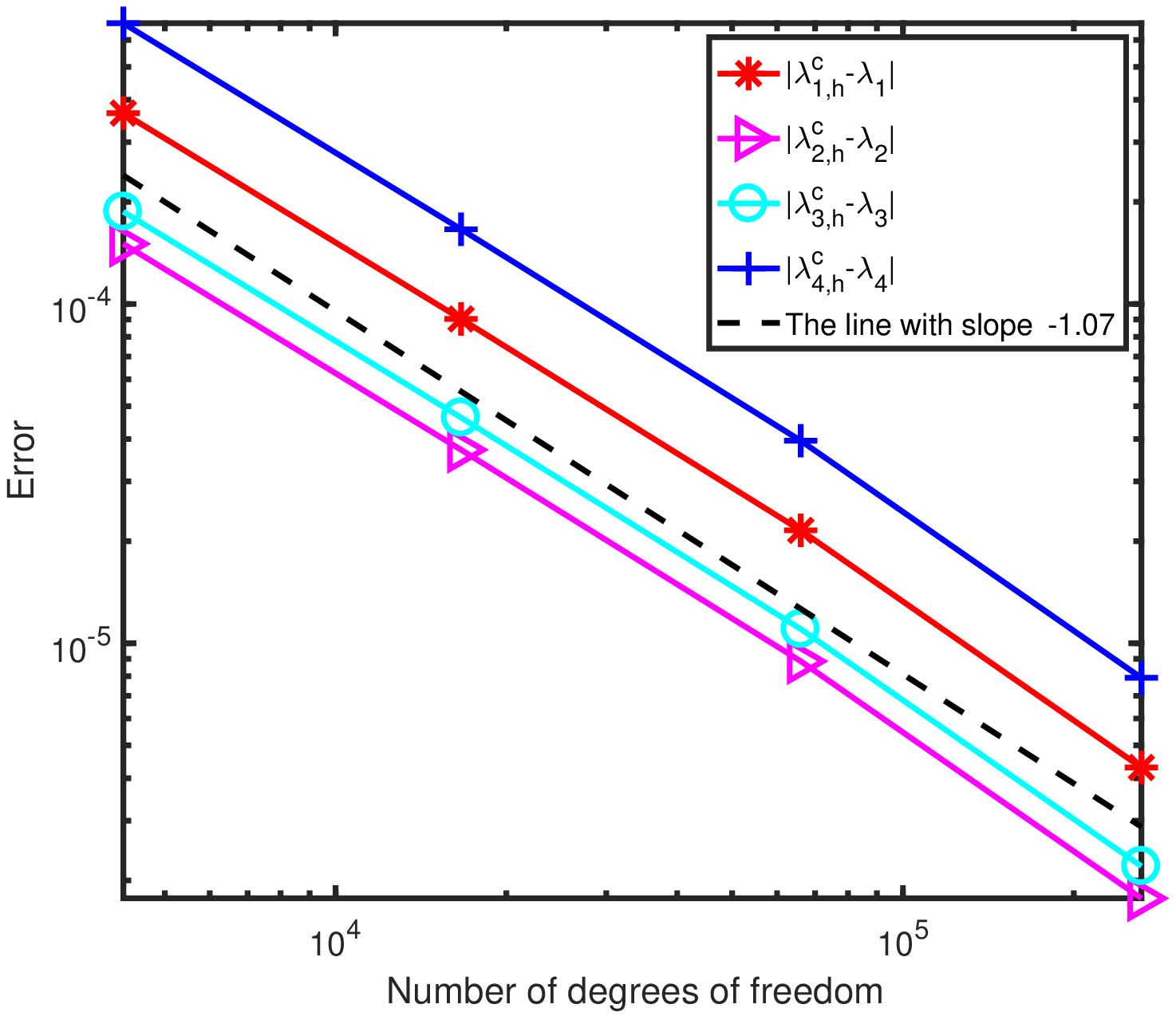}
\includegraphics[width=6.5cm,height=5cm]{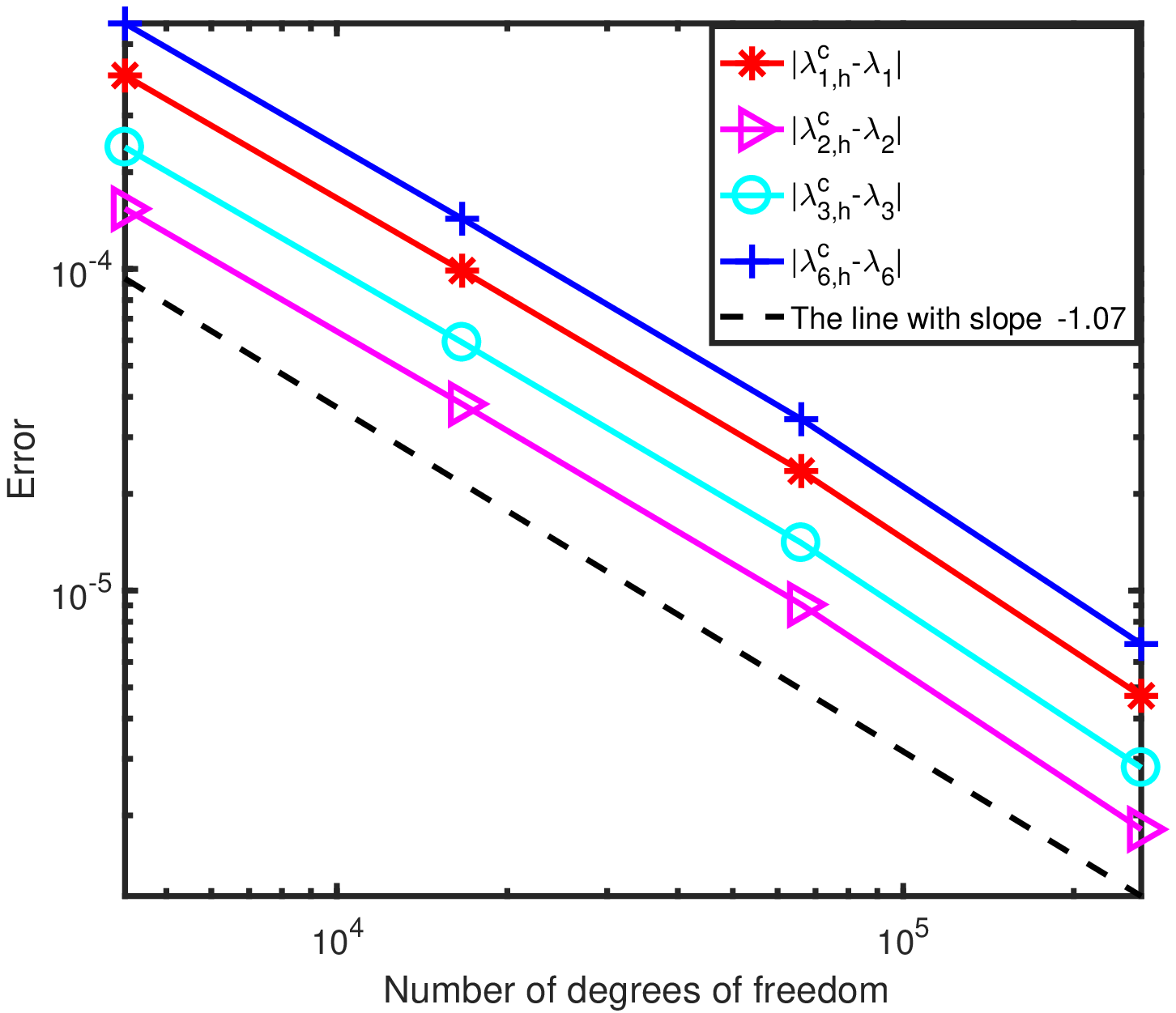}
\caption{\emph{Error curves of eigenvalue approximations on the square(left: $n(x)=4$, right: $n(x)=4+4i$)}}
 \end{figure}
 \begin{figure}[h!]
 \centering
\includegraphics[width=6.5cm,height=5cm]{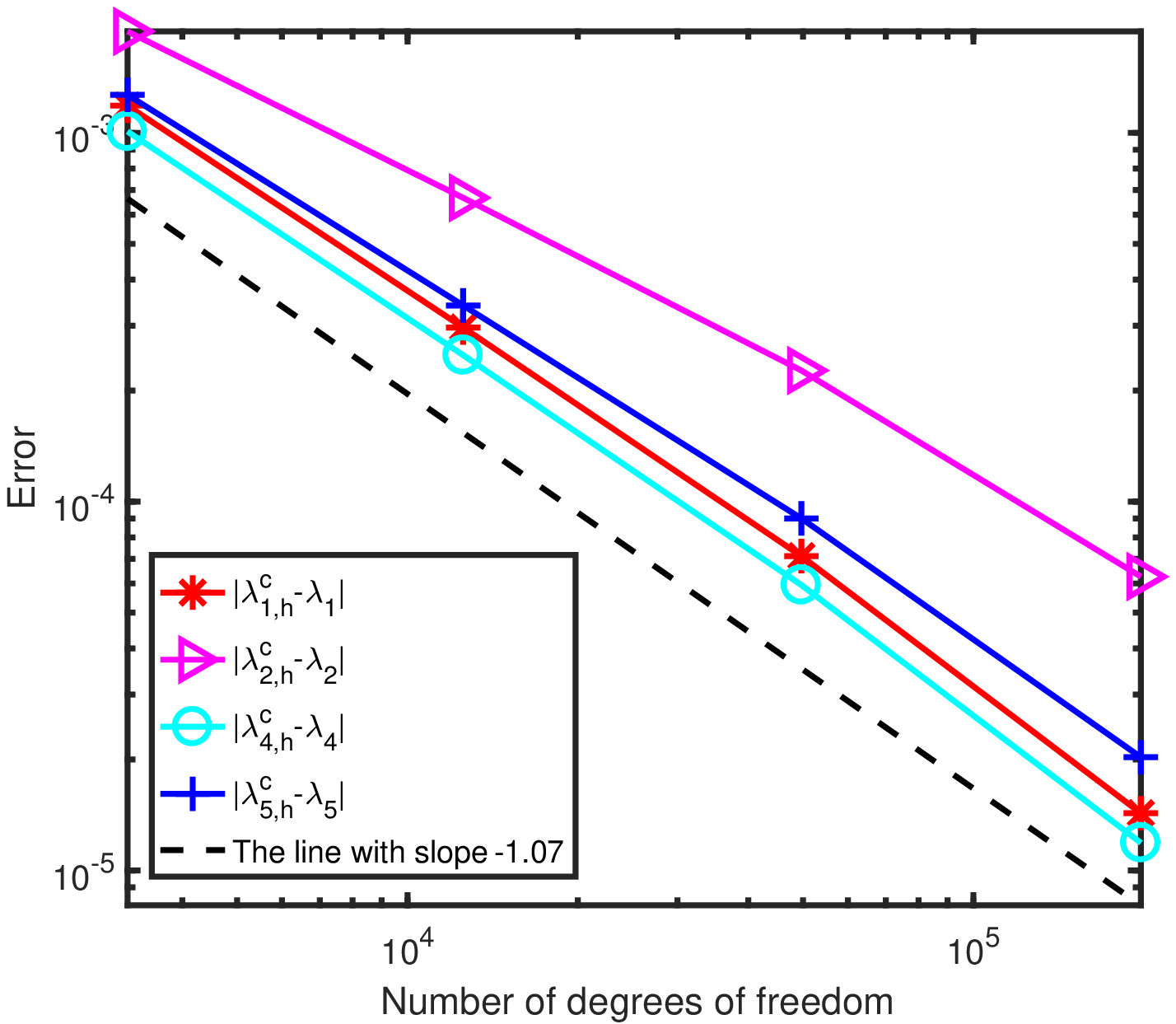}
\includegraphics[width=6.5cm,height=5cm]{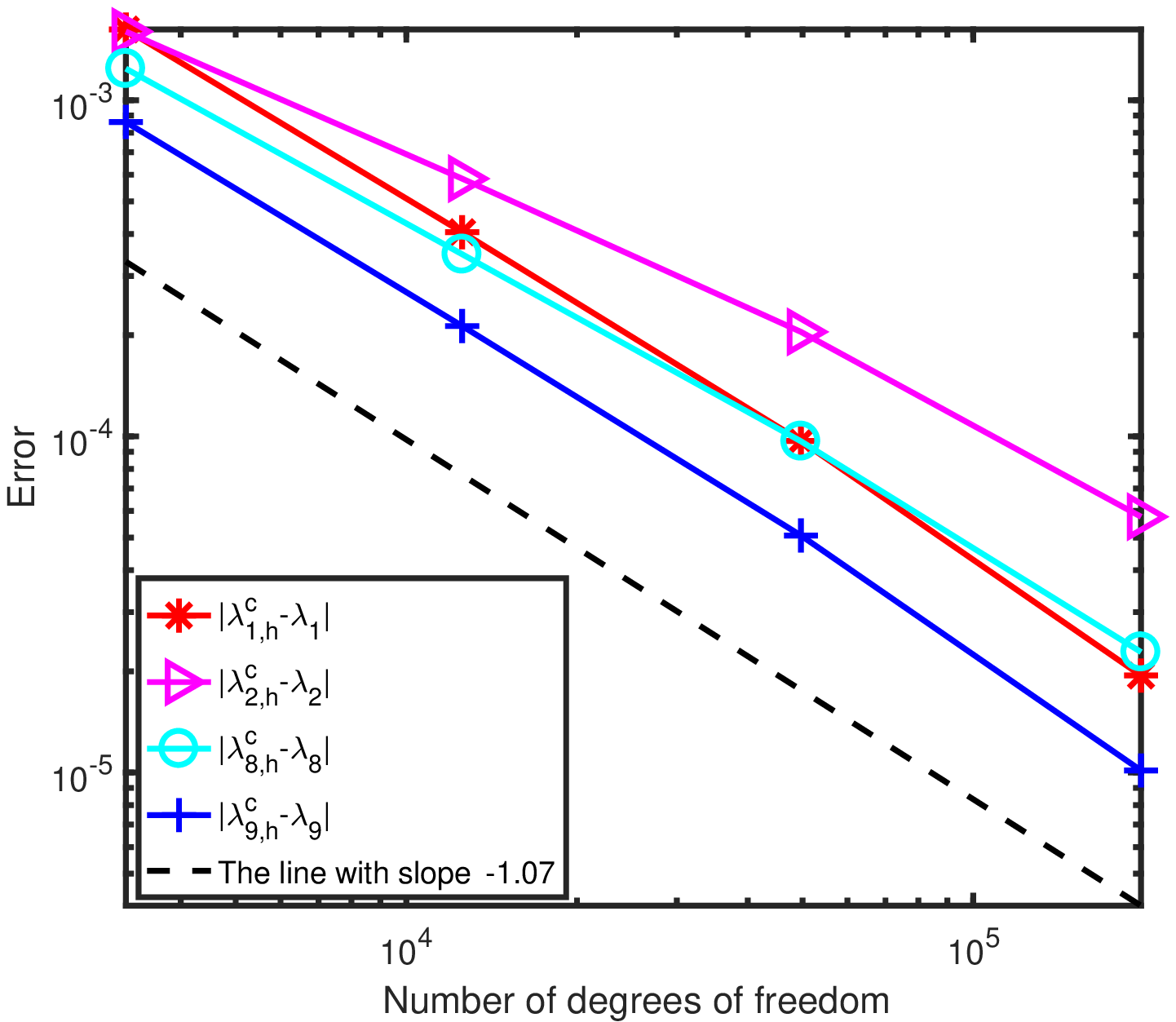}
\caption{\emph{Error curves of eigenvalue approximations on the L-shaped domain (left: $n(x)=4$, right: $n(x)=4+4i$)}}
 \end{figure}
  \begin{figure}[h!]
 \centering
\includegraphics[width=6.5cm,height=5cm]{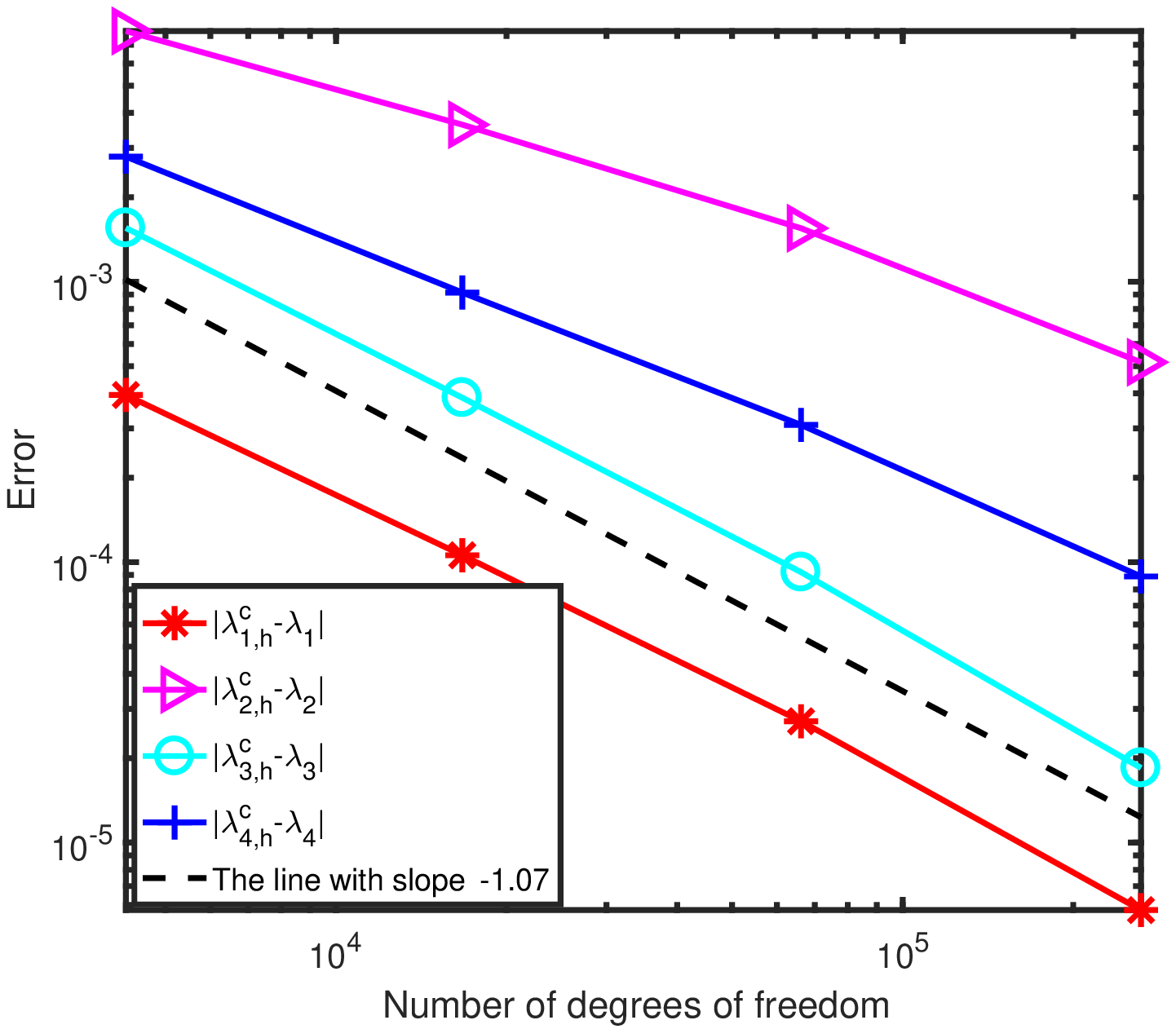}
\includegraphics[width=6.5cm,height=5cm]{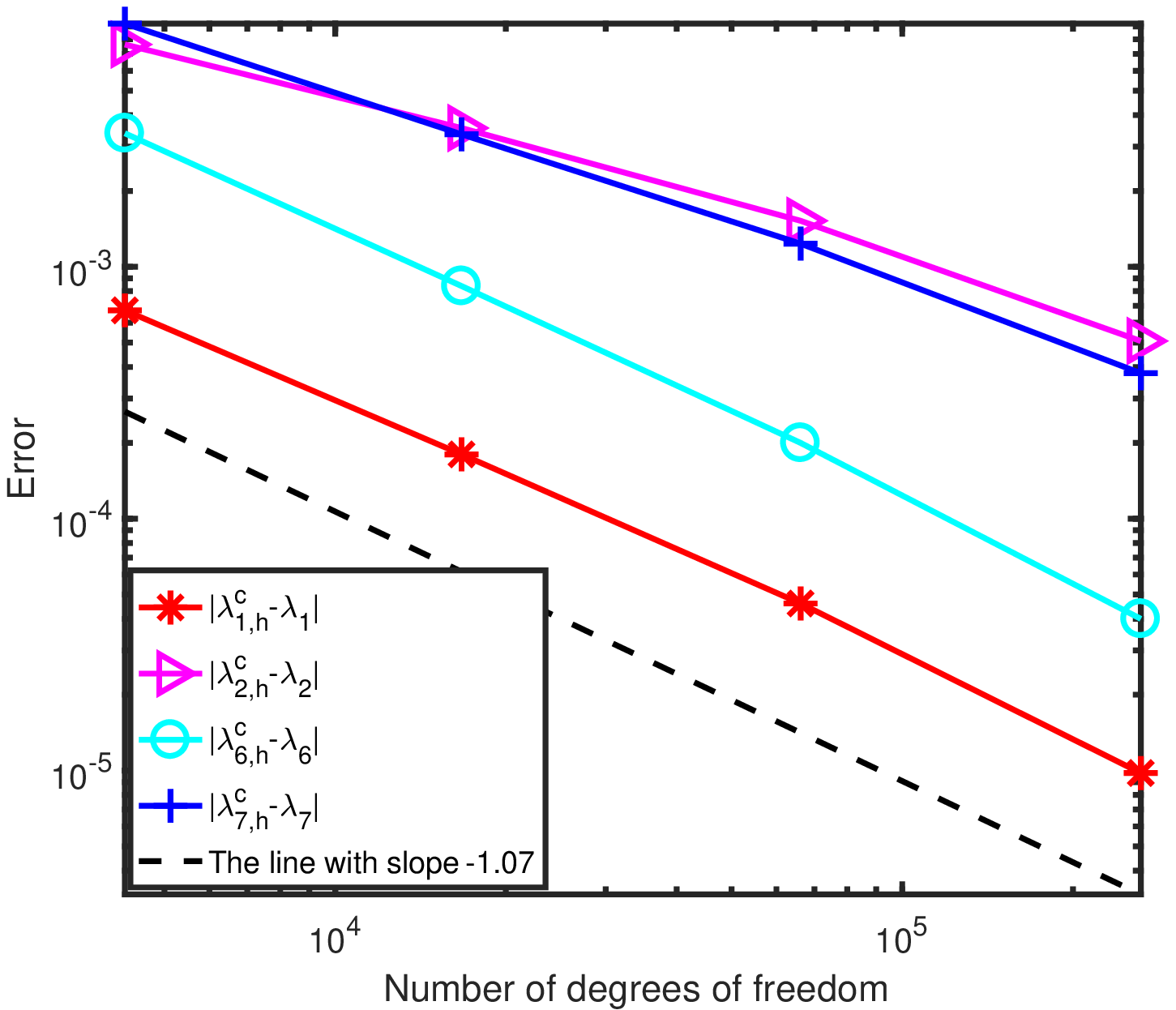}
\caption{\emph{Error curves of eigenvalue approximations on the square with a slit (left: $n(x)=4$, right: $n(x)=4+4i$)}}
 \end{figure}

 \indent Figs 4-6 provide the summations of the errors for the four eigenvalues obtained by algorithm 4.1 and the direct method on each domain. On each figure, we can see that the two curves are almost coincident, which means the multigrid correction scheme can obtain the same optimal error estimates as those by the direct method for the eigenvalue approximations.
\begin{figure}[h!]
 \centering
\includegraphics[width=6.5cm,height=5cm]{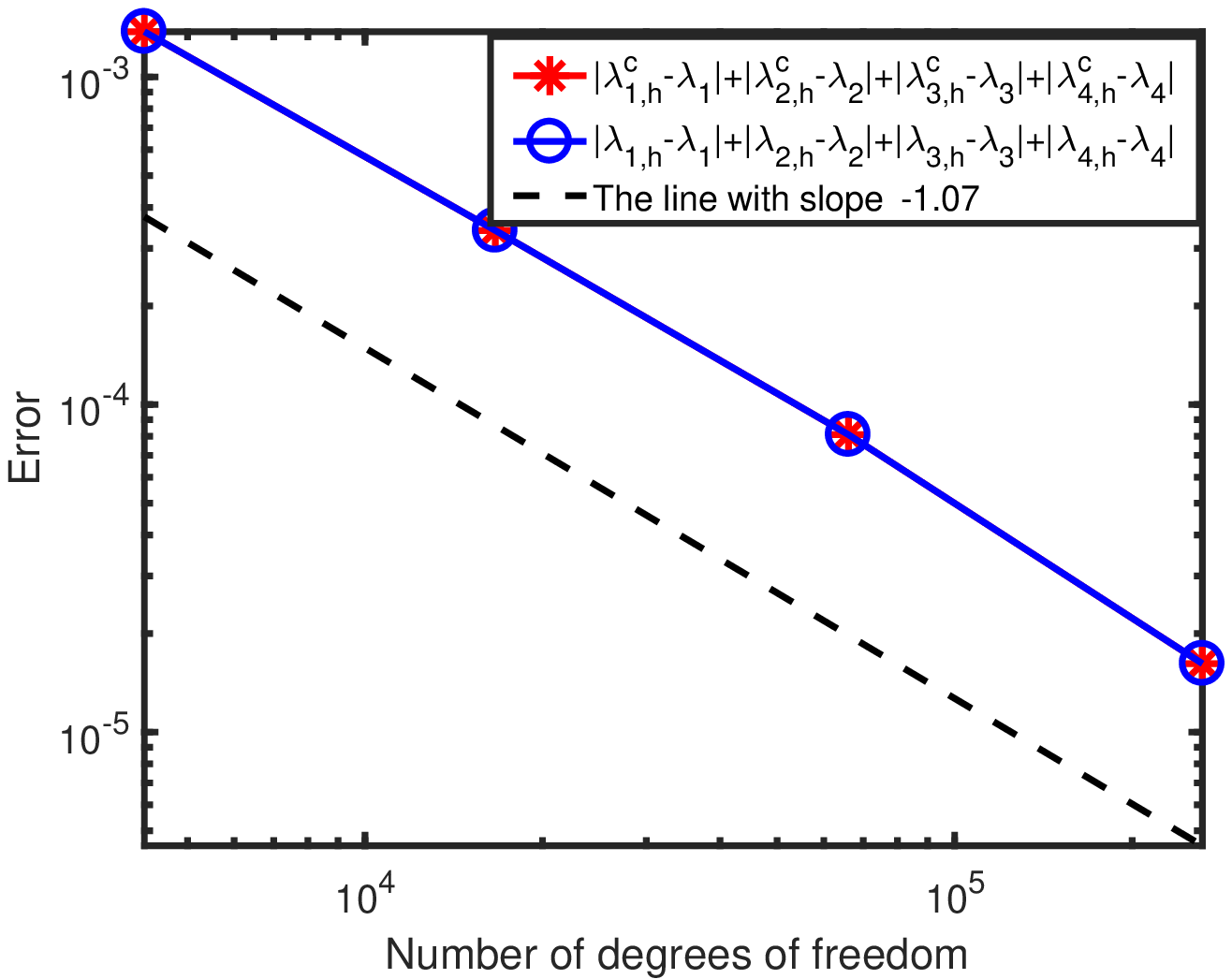}
\includegraphics[width=6.5cm,height=5cm]{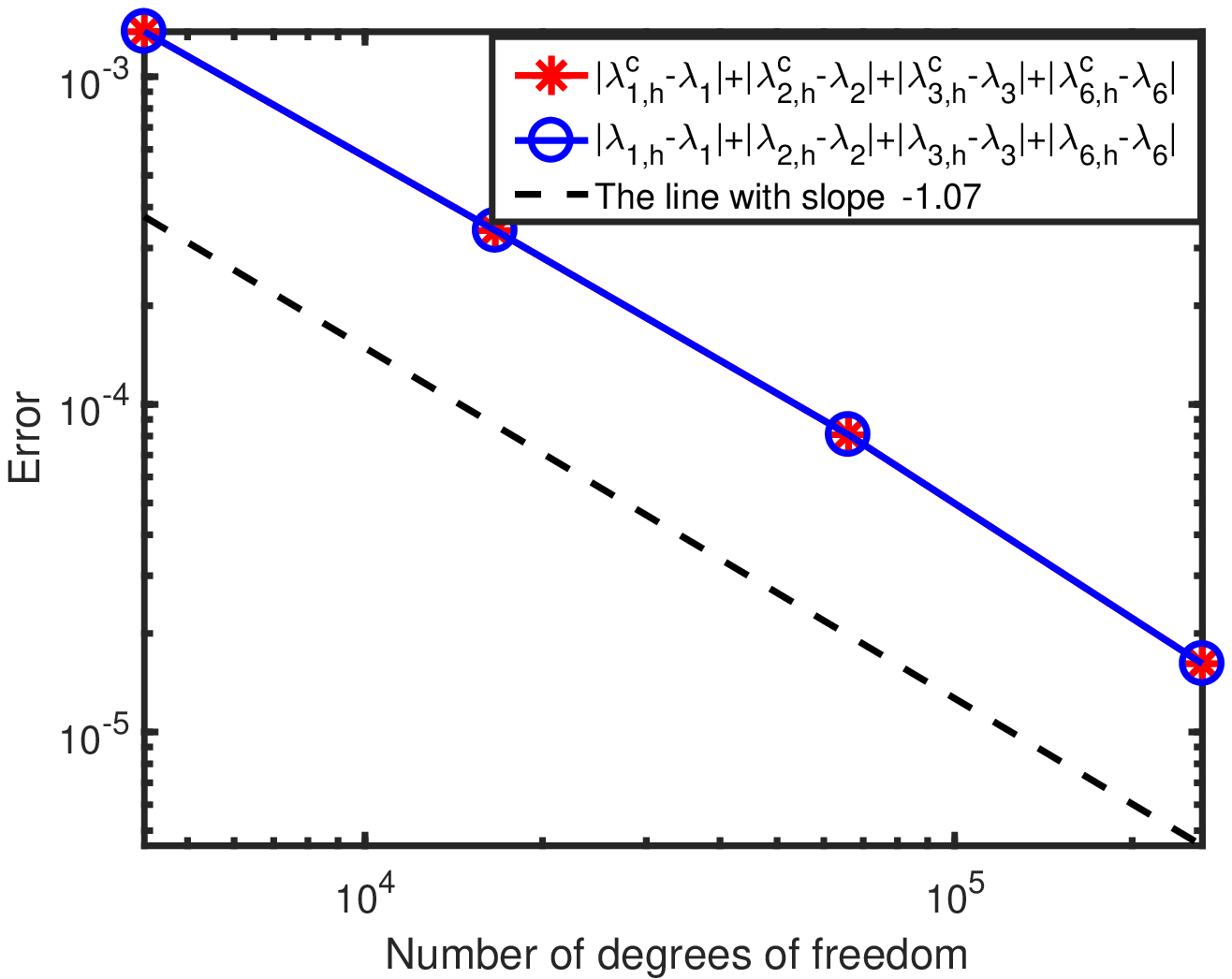}
\caption{\emph{The comparison of the summation of errors for the eigenvalue approximations between the multigrid scheme and the direct method on the square (left: $n(x)=4$, right: $n(x)=4+4i$)}}
 \end{figure}
  \begin{figure}[h!]
 \centering
\includegraphics[width=6.5cm,height=5cm]{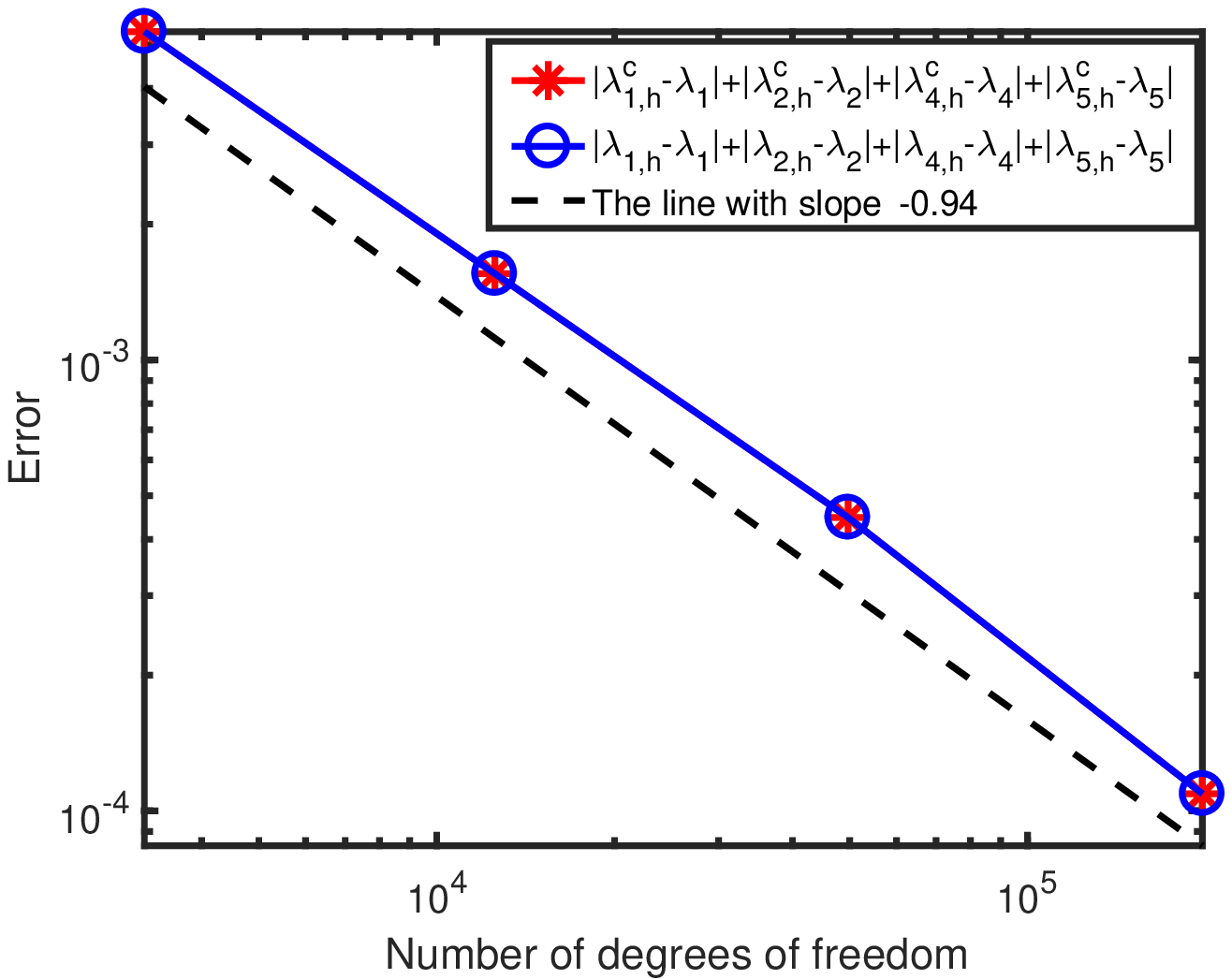}
\includegraphics[width=6.5cm,height=5cm]{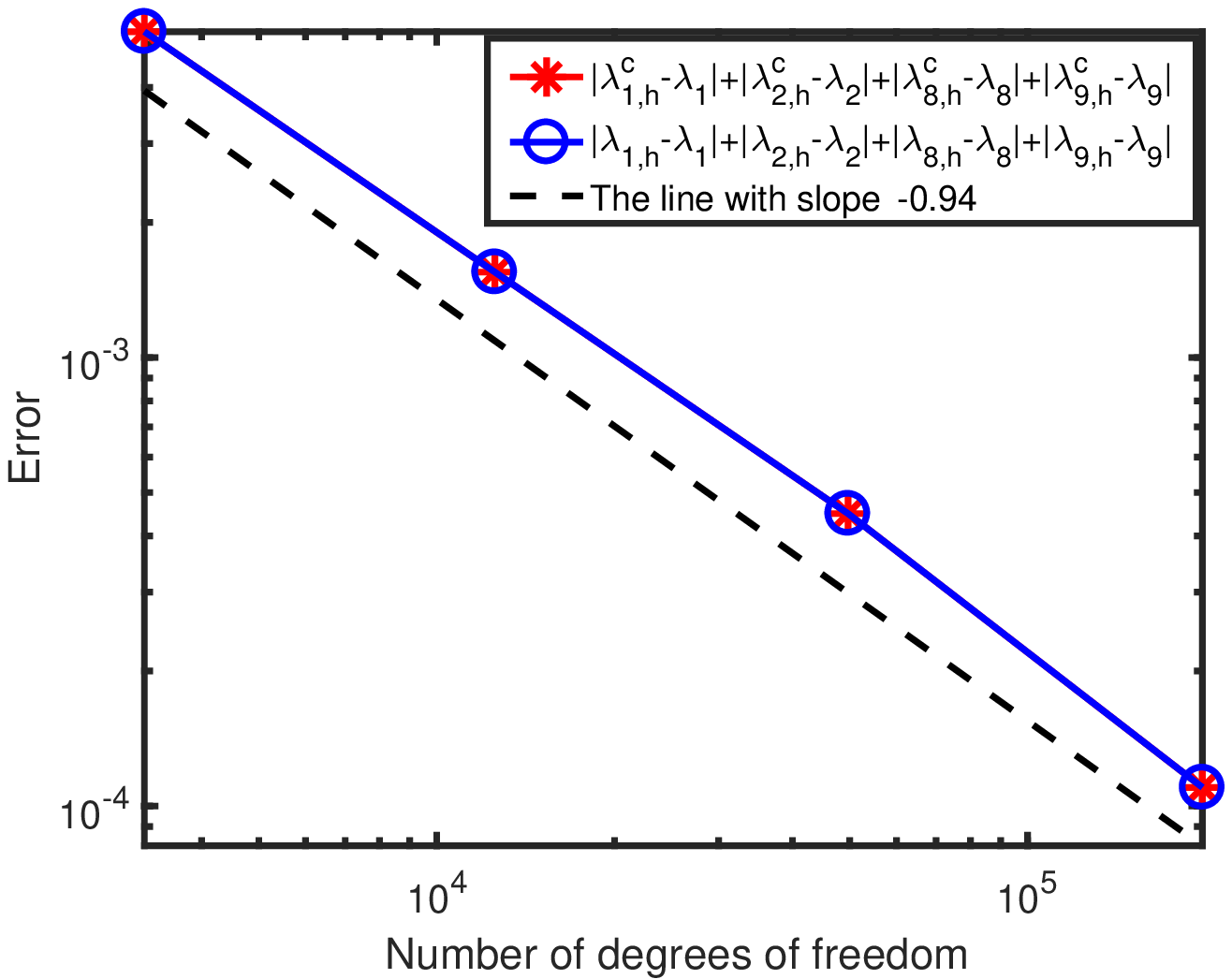}
\caption{\emph{The comparison of the summation of errors for the eigenvalue approximations between the multigrid scheme and the direct method on the L-shaped domain (left: $n(x)=4$, right: $n(x)=4+4i$)}}
 \end{figure}
 \begin{figure}[h!]
 \centering
\includegraphics[width=6.5cm,height=5cm]{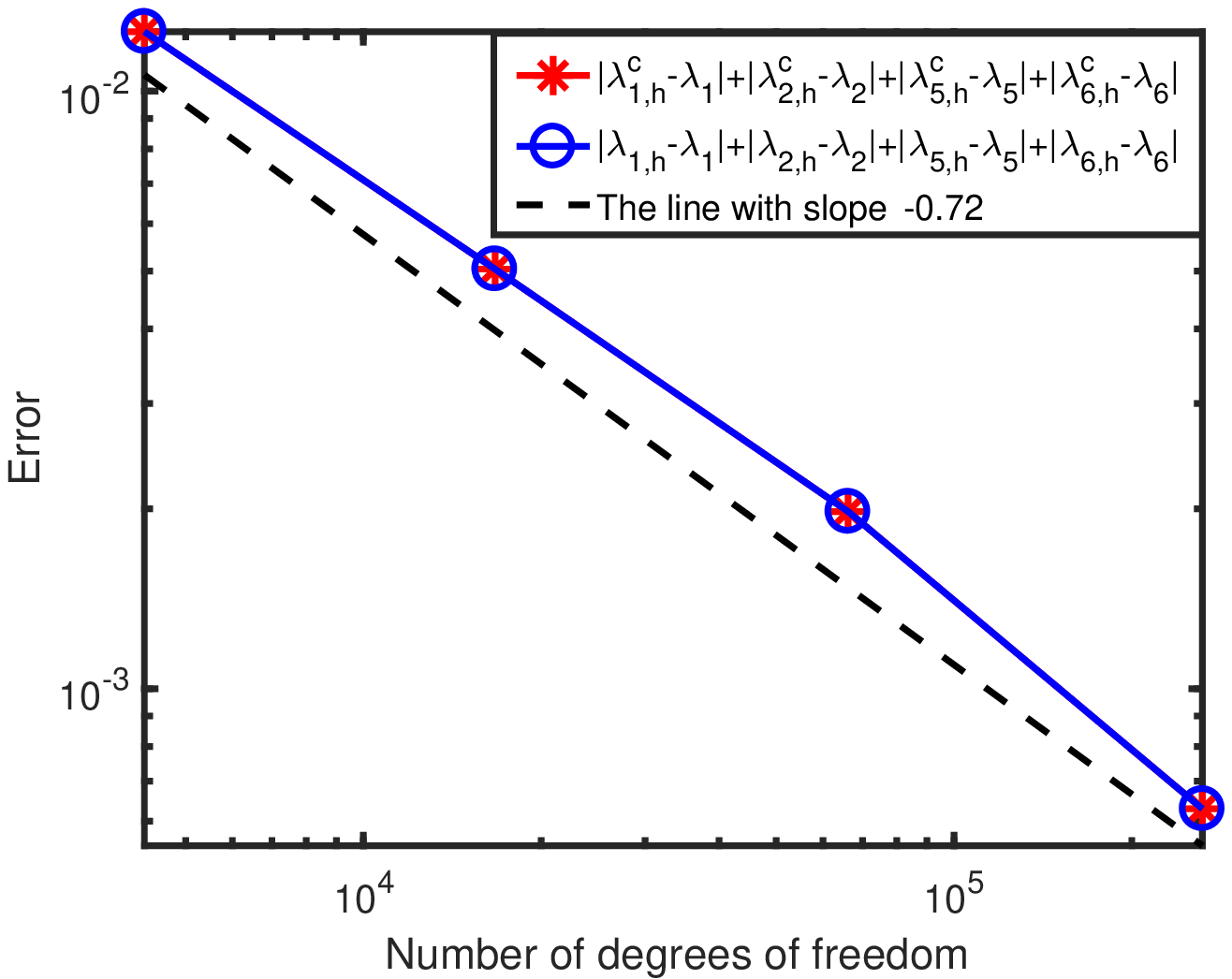}
\includegraphics[width=6.5cm,height=5cm]{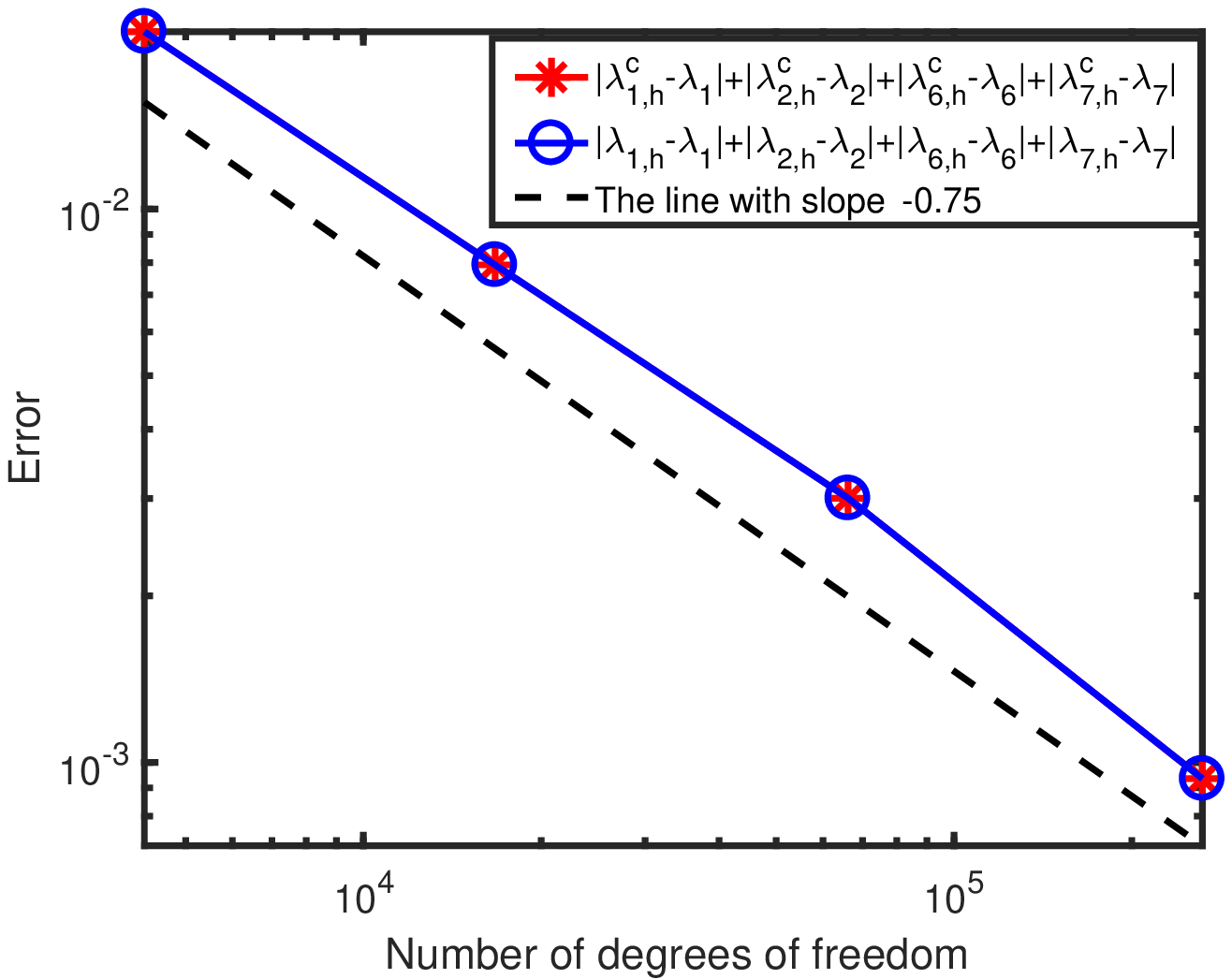}
\caption{\emph{The comparison of the summation of errors for the eigenvalue approximations between the multigrid scheme and the direct method on the square with a slit (left: $n(x)=4$, right: $n(x)=4+4i$)}}
 \end{figure}

  \indent We list the eigenvalue approximations obtained by the multigrid scheme and the direct method in Tables 1-3 when $n(x)=4$ and in Tables 4-6 when $n(x)=4+4i$. From Tables 1-3, we see that there are the same results by the two methods when $n(x)=4$. From Tables 4-6, we see that, when $h$ is the smallest on each table, the eigenvalue approximations can't be computed by the direct method since the computer runs out of memory. With the limited computer memory, our scheme is significant and highly efficient.
  \begin{table}[h!]
  \centering
  \caption{The eigenvalue approximations of (2.3) obtained by Algorithm 4.1 and direct method (square: $n(x)=4$).}
    \begin{tabular}{ccccc}
      \hline\noalign{\smallskip}
    $h$   & $\lambda^c_{1,h}$  & $\lambda^c_{2,h}$& $\lambda^c_{3,h}$&$\lambda^c_{4,h}$\\
    $\frac{2}{512}$ & 2.20250138679  & -0.21225453108  & -0.21225510721  & -0.90806663225  \\
    $\frac{2}{1024}$ & 2.20250569143  & -0.21225275994  & -0.21225290397  & -0.90805872239  \\
   $h$   & $\lambda_{1,h}$  & $\lambda_{2,h}$& $\lambda_{3,h}$&$\lambda_{4,h}$\\
    $\frac{2}{512}$ & 2.20250138680  & -0.21225453108  & -0.21225510721  & -0.90806663225  \\
    $\frac{2}{1024}$& 2.20250569144  & -0.21225275992  & -0.21225290395  & -0.90805872238  \\
   \hline
   \end{tabular}%
\end{table}%
\begin{table}[h!]
  \centering
  \caption{The eigenvalue approximations of (2.3) by Algorithm 4.1 and direct method (L-shaped domain: $n(x)=4$).}
    \begin{tabular}{ccccc}
      \hline\noalign{\smallskip}
   $h$   & $\lambda^c_{1,h}$  & $\lambda^c_{2,h}$& $\lambda^c_{4,h}$&$\lambda^c_{5,h}$\\
    $\frac{2\sqrt{2}}{512}$& 2.53319456612  & 0.85768686246  & -1.08531466335  & -1.09122758504  \\
    $\frac{2\sqrt{2}}{1024}$ & 2.53320886492  & 0.85774947865  & -1.08530278002  & -1.09120730930  \\
     $h$   & $\lambda_{1,h}$  & $\lambda_{2,h}$& $\lambda_{4,h}$&$\lambda_{5,h}$\\
    $\frac{2\sqrt{2}}{512}$ & 2.53319456614  & 0.85768686308  & -1.08531466335  & -1.09122758486  \\
    $\frac{2\sqrt{2}}{1024}$ & 2.53320886479  & 0.85774947872  & -1.08530278003  & -1.09120730928  \\
   \hline
   \end{tabular}%
\end{table}%
\begin{table}[h!]
  \centering
  \caption{The eigenvalue approximations of  (2.3) obtained by Algorithm 4.1 and direct method (square with a slit : $n(x)=4$).}
    \begin{tabular}{ccccc}
      \hline\noalign{\smallskip}
    $h$   & $\lambda^c_{1,h}$  & $\lambda^c_{2,h}$& $\lambda^c_{5,h}$&$\lambda^c_{6,h}$\\
    $\frac{2}{512}$ & 1.48470424178  & 0.46069878294  & -1.89989614930  & -1.92887274492  \\
    $\frac{2}{1024}$ & 1.48470998965  & 0.46121500815  & -1.89987768376  & -1.92878382991  \\
   $h$   & $\lambda_{1,h}$  & $\lambda_{2,h}$& $\lambda_{5,h}$&$\lambda_{6,h}$\\
    $\frac{2}{512}$ & 1.48470424180  & 0.46069878359  & -1.89989614929  & -1.92887274318  \\
    $\frac{2}{1024}$ & 1.48470998967  & 0.46121500835  & -1.89987768376  & -1.92878382943  \\
   \hline
   \end{tabular}%
\end{table}%
 \begin{table}[h!]
  \centering
  \caption{The eigenvalue approximations of (2.3) obtained by Algorithm 4.1 and direct method (square: $n(x)=4+4i$).}
    \begin{tabular}{ccccc}
      \hline\noalign{\smallskip}
   $h$   & $\lambda^c_{1,h}$  & $\lambda^c_{2,h}$& $\lambda^c_{3,h}$&$\lambda^c_{6,h}$\\
     $\frac{2}{512}$& 0.6865580791  & -0.3430478705  & -0.3430446279  & -0.9501192972  \\
          & +2.49529459i & +0.85074449i & +0.85074328i & +0.54009581i \\
    $\frac{2}{1024}$ & 0.6865533933  & -0.3430468763  & -0.3430460656  & -0.9501125093  \\
          & +2.49529414i & +0.850746i & +0.8507457i & +0.54009649i \\
    $h$   & $\lambda_{1,h}$  & $\lambda_{2,h}$& $\lambda_{3,h}$&$\lambda_{6,h}$\\
     $\frac{2}{512}$ & 0.6865580791  & -0.3430478705  & -0.3430446278  & -0.9501192972  \\
          & +2.49529459i & +0.850744489i & +0.8507432795i & +0.540095814i \\
     $\frac{2}{1024}$&  --     &  --     &  --     & -- \\
   \hline
   \end{tabular}%
\end{table}%
\begin{table}[h!]
  \centering
  \caption{The eigenvalue approximations of (2.3) obtained by Algorithm 4.1 and direct method (L-shaped domain : $n(x)=4+4i$).}
    \begin{tabular}{ccccc}
      \hline\noalign{\smallskip}
   $h$   & $\lambda^c_{1,h}$  & $\lambda^c_{2,h}$& $\lambda^c_{8,h}$&$\lambda^c_{9,h}$\\
    $\frac{2\sqrt{2}}{512} $& 0.5143105650  & 0.3969716242  & -1.1594164942  & -1.1423443060  \\
          & +2.88233395i & +1.45891081i & +0.53552365i & +0.52981229i \\
     $\frac{2\sqrt{2}}{1024}$ & 0.5142928934  & 0.3970089008  & -1.1593938106  & -1.1423342849  \\
          & +2.88232587i & +1.45895479i & +0.53552117i & +0.52981075i \\
$h$   & $\lambda_{1,h}$  & $\lambda_{2,h}$& $\lambda_{8,h}$&$\lambda_{9,h}$\\
    $\frac{2\sqrt{2}}{512}$ &  0.5143105650 &  0.3969716223 &  -1.1594164940 &  -1.1423443060  \\
          &  +2.88233395i&  +1.45891081i&  +0.53552365i&  +0.52981229i \\
    $\frac{2\sqrt{2}}{1024}$ &  --     &  --     &  --     & -- \\
   \hline
   \end{tabular}%
\end{table}%
\begin{table}[h!]
  \centering
  \caption{The eigenvalue approximations of (2.3) obtained by Algorithm 4.1 and direct method (square with a slit : $n(x)=4+4i$).}
    \begin{tabular}{ccccc}
      \hline\noalign{\smallskip}
   $h$   & $\lambda^c_{1,h}$  & $\lambda^c_{2,h}$& $\lambda^c_{6,h}$&$\lambda^c_{7,h}$\\
    $\frac{2}{512}$ & 0.919316438  & 0.291737235  & -2.859202716  & -2.850238759  \\
          & +1.77078218i & +0.99939462i & +0.50477279i & +0.49299969i \\
     $\frac{2}{1024}$ & 0.919307780  & 0.292183070  & -2.859162853  & -2.849868193  \\
          & +1.77078671i & +0.9996367i & +0.5047693i & +0.4930741i \\
    $h$   & $\lambda_{1,h}$  & $\lambda_{2,h}$& $\lambda_{6,h}$&$\lambda_{7,h}$\\
    $\frac{2}{512}$ &  0.919316438 &  0.291737232 &  -2.859202716 &  -2.850238737  \\
          &  +1.77078218i&  +0.99939462i&  +0.50477279i&  +0.49299969i\\
    $\frac{2}{1024}$ &  --     &  --     &  --     & -- \\
   \hline
   \end{tabular}%
\end{table}%
\section{Conflict of Interests}
\indent The authors declare that there is no conflict of interests
regarding the publication of this paper.
\section{Acknowledgments}
\indent This work is supported by the National Science Foundation
of China (Grant no. 11761022) and the Thousands of levels of innovative talent support(2015).

\end{document}